\documentclass[12pt,mathpaper]{article}
\usepackage[margin=1in]{geometry}
\usepackage{graphicx}
\usepackage{bbding}
\usepackage{perpage}
\usepackage{mathtools}
\usepackage{multicol}
\usepackage{calc}
\usepackage{amsmath, ulem}
\usepackage{amsxtra}
\usepackage{enumerate}
\usepackage{tikz-cd}
\usepackage[shortlabels]{enumitem}
\usepackage[thicklines]{cancel}
\usepackage{amssymb,wasysym}
\usepackage{amsfonts}
\usepackage{pifont}

\usepackage{amsthm}
\usepackage[sans]{dsfont}
\usepackage{xcolor,cancel}
\usepackage{verbatim}

 \definecolor{lightgrey}{rgb}{0.7,0.7,0.7}
\usepackage{textcomp}
\usepackage{mathrsfs}
\usepackage[colorlinks]{hyperref}

\newtheorem{thm}{Theorem}[section] 
\newtheorem{lem}[thm]{Lemma}

\newtheorem{prop}[thm]{Proposition} 
\newtheorem{assu}[thm]{Assumption} 
\newtheorem{rmk}[thm]{Remark}

\newcommand{\Diff}{\mathsf{Diff}(D)}
\newcommand{\MPD}{\mathsf{SDiff}(D)}
\newcommand{\Dens}{\mathsf{Dens}(D)}

\newcommand{\sG}{\mathsf{G}}
\newcommand{\bId}{\mathbf{Id}}
\newcommand{\leps}{\bigg\lvert_{\varepsilon=0}}
\newcommand{\Hess}{\mathsf{Hess}}

\newcommand{\bP}{\mathbf{P}}
\newcommand{\grad}{\mathsf{grad}}

\newcommand{\Ker}{\mathsf{Ker}}
\newcommand{\Ran}{\mathsf{Ran}}
\newcommand{\Tr}{\mathsf{Tr}}
\newcommand{\bL}{\mathbf{L}}
\newcommand{\bmat}{\mathbf{M}}

\usepackage{caption}
\title{A Lagrangian Approach to the Inhomogeneous Incompressible Euler Equation}
\author{Anping Pan}

\begin{document}

\maketitle

\begin{abstract}
 In this paper, we study the  inhomogeneous  incompressible Euler equation (IIE in short) from a Lagrangian perspective. We establish a geodesic description of this equation and discuss the associated geometric structures. We also find the derivation of IIE from the Hamilton-Pontryagin action principle and derive the corresponding Lagrangian formulation. A byproduct is a new vorticity formulation of IIE. We also prove the Lagrangian analyticity of IIE using our Lagrangian representation formula.  
\end{abstract}

\tableofcontents

\section{Introduction}

The geometry of hydrodynamics has been a vivid research direction since the seminal work of Arnold \cite{Ar66} and Ebin-Marsden \cite{EM70}, where the motion of perfect fluids in a smooth bounded domain $D$ is interpreted as geodesic in the infinite dimensional Riemannian manifold $\MPD$ of volume-preserving diffeomorphisms of $D$, which is equipped with the right-invariant Riemannian metric given by the kinetic energy of the fluid. This geometric perspective has been extensively generalized since then; we refer to \cite{AK98}\cite{HMR98}\cite{MR99} for further development of this topic, and \cite{KMM20} for a recent survey on geometric hydrodynamics including compressible fluids.

From a mechanical viewpoint, Arnold's geodesic formulation is a natural infinite-dimensional generalization of the classical least action principle in Lagrangian mechanics. Flipping to the Hamiltonian side, we have an analogous Hamilton-Pontryagin type variational principle for hydrodynamic models (see \cite{BCHM00}\cite{Ho21}), which leads to the Lagrangian-Eulerian formulation of fluid equations \cite{Co01a}. We remark that the above perspective admits a natural generalization which effectively applies to viscous hydrodynamic models. We refer to \cite{CI08} for the Navier-Stokes case and to  \cite{MP25} for an extended discussion. However, all the above discussions are restricted to fluid models with mechanical Lagrangians, where  velocity and position dependence are decoupled and appear only in kinetic and potential energy respectively.

A natural example which falls outside the class of mechanical Lagrangians is the following inhomogeneous version of the incompressible Euler equation (IIE in short), which is an extensively studied model in past decades:
\begin{equation}\label{IIE}
\left\{
\begin{aligned}
&\partial_t(\rho u)+\nabla\cdot(\rho u\otimes u)+\nabla p=0\\
&\nabla\cdot u=0\\
&\partial_t\rho+u\cdot\nabla\rho=0
\end{aligned}
\right.
\end{equation}
The above equation describes the motion of incompressible fluids with variable density.  Here, the inhomogeneity of the fluid forces the kinetic energy density to depend on both position and velocity variables:
\begin{equation*}
\frac{1}{2}\lvert u\rvert^2\rho=\frac{1}{2}\lvert u\rvert^2 \rho_0\circ X^{-1},\quad \dot X=u(X).   
\end{equation*}

The above equation was systematically studied in the celebrated monograph \cite{[Lio96]}. Analysis of the equation has been done in various functional settings in a series of recent works, see \cite{Da06}\cite{[Da10]}\cite{[DF11]}. The main focus of the works listed above is to address the well-posedness issues in various function spaces and to analyze IIE from an Eulerian viewpoint using Harmonic analysis tools. However, contrary to the well-known Yudovich theorem of $2$-$d$ homogeneous Euler's equation, the global well-posedness issue of IIE in $2$d is still open. In recent years, there have also been works devoted to Beale–Kato–Majda-type (cf.\cite{[BKM84]} ) blow-up/continuation criteria for IIE, we refer the readers to \cite{[BLS20]}\cite{Ch03}\cite{[Fa25]}.

However, different from the analytical side, the geometric structure and Lagrangian perspective of the inhomogeneous Euler equation is much less studied. To the best of the author's knowledge, the only relevant literature in this direction is \cite{[Mar76]}, which studies the geodesic formulation of IIE on $\MPD$, as a generalization of the celebrated geodesic formulation of the incompressible Euler equation. The Hamiltonian aspects and Lagrangian formulation (a Weber-type formula, see \cite{Co01a}) of IIE, to the best of the author's knowledge, are absent in the current literature. 

To fill the gap in the literature, in this paper we establish a Lagrangian approach to the analysis of IIE, which stems from geometric consideration. We discuss the Hamiltonian formulation of the equation and derive a Lagrangian-Eulerian representation of the solution of IIE. The advantages of such a Lagrangian perspective are twofold: First, the Lagrangian formulation motivates us to find a vorticity formulation of IIE, which totally eliminates the pressure; Second, benefiting from the Lagrangian perspective, we obtain an explicit recursive formula for the time-series coefficients of the Lagrangian displacement, and demonstrate Lagrangian analyticity of the solution accordingly. 

We believe that the vorticity formulation we obtained for IIE may uncover new approaches to the study of IIE. The key takeaway here is that: Geometrically, vorticity in incompressible fluid models is naturally given by the curl of the \textit{momentum} rather than the velocity, due to the dual Lie transport (coadjoint motion) nature of momentum, together with the following commutation relation, which encodes the vorticity transport:
\begin{equation*}
 \nabla\times\mathscr L_u^*=\mathscr L_u \nabla\times,\quad \mathscr L_u^*v:=u\cdot\nabla v+\nabla^* uv.   
\end{equation*}
Thus, there is a natural difference between the vorticity of homogeneous and inhomogeneous quantities, based on geometric considerations. As the vorticity formulation plays an important role in the study of various problems of incompressible Euler equation, for instance the finite time blow up, continuation criteria and stability of steady states, the author believes that our vorticity formulation may be helpful in investigating analogous problems for IIE. In our companion paper \cite{Pa26+}, based on our vorticity formulation, we study the lifespan and continuation criteria of IIE via a Lagrangian approach.

The structure of this paper is organized as follows: In section 2, we study the geodesic formulation of IIE with respect to a weighted Riemannian metric on $\MPD$, where we interpret the weighted Leray projector as the second-fundamental form. The sectional curvature and Lagrangian stability are also discussed. In section 3, we derive IIE from a Hamilton-Pontryagin type variational principle, and obtain a Lagrangian-Eulerian formulation of IIE as a Lagrangian formula for the solution of the Hamiltonian ODE given by the variational principle. We also derive the vorticity formulation from our geometric construction, and recover the local existence results in our framework. Finally, in section 4, using our Lagrangian representation, we prove the time-analyticity of Lagrangian trajectories of the IIE.

\section*{Acknowledgements} The author thanks V.Vicol for stimulating discussions which initiate the study of Lagrangian analyticity of IIE. The author thanks T.Drivas for pointing out the ongoing work \cite{DG26+} and for helpful discussions. The author thanks B.Khesin, H.N.Lopes, M.Lopes, A.Mazzucato and A.Novikov for stimulating discussions and valuable suggestions. This work was partially supported by the US National Science Foundation grants DMS-1909103 and DMS-2206453.

\section{Geometry of \texorpdfstring{$\MPD$}{MPD} with Weighted Metric}

In this section, we study the Riemannian structure of the volume-preserving diffeomorphism group, which is equipped with a weighted Riemannian metric and naturally fits the inhomogeneous incompressible Euler equation. Our presentation in this section mainly follows and generalizes the discussions in the unpublished notes \cite{DrWS25}, which generalizes constrained geodesic, second fundamental form and sectional curvature of level set submanifold from finite to infinite dimension.

\subsection{Some Conventions}

\qquad We fix some frequently used notations. Let $D$ be the physical domain of our fluids. Unless otherwise specified, we assume $D=\mathbb T^d$, that is the unit box $[0,1]^d$ equipped with periodic boundary condition. We remark that most our results also hold in the setting of smooth compact Riemannian manifold or $\mathbb R^d$, and works under $H^{s+1}$ Sobolev regularity with $s>d/2$, after suitable technical modification (e.g. replacing mean-zero condition with slip boundary condition for compact manifold $D$, or with decay condition for $D=\mathbb R^d$). 

Throughout this paper, we denote by $\Diff$($\MPD$) the Lie group of (volume-preserving) diffeomorphisms in $D$. We denote by $\mathfrak X(D)$, $\mathfrak X_0(D)$, $\mathfrak X_\sigma(D)$ the space of smooth vector fields, mean-free vector fields and divergence-free  vector fields on $D$. 

In this section, we will focus on the geometric formulation of IIE and only consider the smooth case on torus for simplicity. Later, for analytical results like well-posedness and analyticity, we will work in Sobolev regularity class.  We will use the notation $H_0^k(D)$ for space of mean-zero $H^k(D)$ functions, and we will frequently use the notation $H^s$ without specifying the domain and range of the function. We denote by $\mathbf{P}$ the standard Leray projector on $\mathfrak X(D)$, which projects vector fields to its div-free part. In compact form:
\begin{equation*}
  \bP v=v-\nabla\Delta^{-1}\nabla\cdot v.  
\end{equation*}

It is well-known that $\bP$ is a zero-order multiplier that preserves the mean value, Hence, it is bounded on each Sobolev space $H_0^k(D)$ and $\bP: \mathfrak X_0(D)\to \mathfrak X_0(D)\cap \mathfrak X_\sigma(D)$.

We will use the letter $X$ for elements in $\MPD$, and letter $A$ for the inverse of $X$, which is the back-to-label map associated with flow map $X$ in fluid mechanics context.

 We denote by $\Dens$ the space of smooth densities in $D$, that is:
 \begin{equation*}
 \Dens=\bigg\{\rho\in C^\infty(D): \int_D \rho(x)dx=1\bigg\} .   
 \end{equation*}
 
 We will fix a reference density $\rho_0\in\Dens$, which is required to satisfy the following assumption throughout this paper:
\begin{assu}[Assumption (A)]\label{Assumption A}
 There exists $C_1>1$ such that
\begin{equation*}
 0<C_1^{-1}\le \rho_0(a)\le C_1<\infty\quad\text{ for all }a\in D .  
\end{equation*}
\end{assu}
The above assumption \ref{Assumption A} is posed to avoid concentration or vacuum in fluid density. We denote by 
\begin{equation*}
 \rho=X_\sharp\rho_0=\rho_0\circ A,\quad X\in\MPD,   
\end{equation*}
which is the dynamical density transported by some rearrangement $X$. We notice that the assumption \ref{Assumption A} is propagated under the action of rearrangements. Hence, also propagated by transport equation. We may sometimes add the subscript $X_t$ to highlight its time dependence.

We will frequently deal with variable coefficient second-order elliptic operators. Fix a $\rho\in\Dens$, we consider the following weighted Laplacian operator and the associated Poisson equation:

\begin{equation*}
 \Delta_\rho \varphi:=\nabla\cdot(\rho^{-1} \nabla \varphi)=f\quad\text{ in }D.   
\end{equation*}

Assumption \ref{Assumption A} guaranties that  $\Delta_\rho$ is an uniformly elliptic operator. Here, we assume the forcing term $f$ satisfies the mean-zero condition:
\begin{equation*}
 \int_D f(x)dx=0.   
\end{equation*}
Then, $\Delta_\rho$ is invertible in the space of mean-zero smooth functions. We will denote by $\mathcal G_\rho$ the Green's function associated with $\Delta_\rho$, and the following holds:
\begin{equation*}
 \int_D \mathcal G_\rho(x,y)dx=0,\quad \mathcal G_\rho(x,y)=\mathcal G_\rho(y,x),\quad  \varphi(x)=\int_D \mathcal G_\rho(x,y)f(y)dy.  
\end{equation*}

We will use the letter $C$ to denote positive constants which may depend on additional parameters and change from line to line. $A\lesssim B$ means there exists $C>0$ which might depend on other harmless parameters such that $A\le CB$.

We end this subsection with the well-known Jacobi-Liouville formula, which will be repeatedly used in the sequel. 
\begin{lem}\label{Abel-Jacobi-Liouville Formula}
Let $\bmat: [0,1]\to \mathsf{GL}(d)$ be a smooth curve in the general linear group $\mathsf{GL}(d)$, then we have:
\begin{equation}\label{Jacobi Identity}
 \frac{d}{dt}\det(\bmat)=\Tr(\bmat^{-1}\dot \bmat)\det(\bmat);   
\end{equation}
and
\begin{equation}\label{derivative of A^-1}
 \frac{d}{dt}\bmat^{-1}= -\bmat^{-1}\dot  \bmat \bmat^{-1}  .
\end{equation}    
\end{lem}
\begin{proof}[Proof of Lemma~\ref{Abel-Jacobi-Liouville Formula}]

We first notice \eqref{derivative of A^-1} is a simple consequence of the chain rule:
\begin{equation*}
0=\frac{d}{dt}(\bmat\bmat^{-1})=\dot \bmat\bmat^{-1}+\bmat\frac{d}{dt} \bmat^{-1}.    
\end{equation*}
Rearranging the above equality, we obtain \eqref{derivative of A^-1}.  

 To show \eqref{Jacobi Identity}, the standard approach is to use minor expansion.  Here, we provide a simple proof suggested by \cite{Kh25}.

 Notice that $\det \bmat$ is the Wronskian of matrix $\bmat$, since the matrix ODE that $\bmat$ satisfies reads
 \begin{equation*}
     \dot\bmat=\bmat (\bmat^{-1}\dot\bmat)=:\bmat \tilde{\bmat},
 \end{equation*}
 we apply the Abel's formula for Wronskian (c.f. \cite{TeschlODE2012}) to conclude:
 \begin{equation*}
  \frac{d}{dt}\det \bmat=\Tr(\tilde{\bmat})\det\bmat_t=\Tr(\bmat^{-1}\dot{\bmat})\det\bmat. 
 \end{equation*}
   
\end{proof}

\subsection{The Riemannian Structure}

Our geometric description relies on the interpretation of $\MPD$ as a level set submanifold of the full diffeomorphism group $\Diff$ with infinite dimension and co-dimension, and the inhomogeneous Euler equation will be derived as geodesic on the submanifold, with respect to some weighted metric. This formalism closely follows \cite{DrWS25}. (See also the ongoing work \cite{DG26+}, dealing with dynamics on rapidly oscillating constraint submanifold in infinite dimensions.)

More precisely, we rewrite $\MPD$ as follows:
\begin{equation*}
\MPD=\bigcap_{a\in D}\{X\in \Diff :f_a(X)=0\},\quad f_a(X)=\langle \det(\nabla X)-1,\delta_a\rangle_*    .
\end{equation*}
Here, the duality pairing $\langle\cdot,\cdot\rangle_*$ is the Riesz pairing between $C_0$ functions and measures.

Now for $X\in\MPD$, to see the tangential direction of the level sets of $f_a$, we apply \eqref{Jacobi Identity} to compute
\begin{equation*}
\frac{d}{d\varepsilon}\leps f_a(X+\varepsilon\delta X)=\langle\Tr((\nabla X)^{-1}\nabla\delta X)\det(\nabla X),\delta_a\rangle_*    
\end{equation*}
\begin{equation*}
=[\nabla\cdot(\delta X\circ X^{-1})](X(a))=0,    
\end{equation*}
which implies that
\begin{equation*}
 T_X\MPD=\{v\circ X: v\in \mathfrak X_\sigma(D)\} .  
\end{equation*}

Define the following weighted metric:
\begin{equation}
\sG_X(\xi,\eta)=\int_D (\xi,\eta)(x)\rho_0(x)dx=\int_D (\xi\circ X^{-1},\eta\circ X^{-1})\rho(x)dx    .
\end{equation}

Under the above metric, the normal space $N_X\MPD=T_X^\perp\MPD$ is given by
\begin{equation}
  N_X\MPD=\{\rho_0^{-1}\nabla \varphi(X):\varphi\in H^1(D)\}  .
\end{equation}
Hence, we can define the following map that generates normal vectors:
\begin{equation}
\mathcal N_{X}:\varphi\in H^1(D)\to -\rho_0^{-1}\nabla \varphi(X)=-\rho_0^{-1}(\nabla^*X)^{-1}\nabla(\varphi\circ X)\in N_X\MPD    .
\end{equation}

We now extend the domain of the above normal generator $\mathcal N$ to distributions. Let $\mathcal M$ be the set of Borel measures on $D$ with finite variation, for any $\mu\in\mathcal M$ we define the following co-normal functional 
\begin{equation*}
\mathcal N_X^*(\mu)\in T_X^*\MPD,\quad \mathcal N_X^*(\mu):= \delta X\in T_X\Diff\to \langle \mu,\nabla\cdot(\delta X\circ A)\rangle_{*}   . 
\end{equation*}
which in distributional sense reads:
\begin{equation*}
\mathcal N_X^*(\mu)=-\nabla \mu(X).
\end{equation*}
Moreover, the annihilating property $T_X\MPD\subset \Ker(\mathcal N_X^*(\mu))$ is satisfied. 

Moreover, in this context, the musical isomorphism (inertia operator in the Language of mechanics) is $\alpha\to \rho_0^{-1}\alpha$. Hence, the vector representative $\mathcal N_X$ of covector $\mathcal N_X^*$ is:
\begin{equation*}
 \mathcal N_X(\mu):=\rho_0^{-1}\mathcal N_X^*(\mu)=-\rho_0^{-1}\nabla \mu(X)=-\rho_0^{-1}(\nabla^* X)^{-1} \nabla(X_\sharp\mu) , 
\end{equation*}
which should be interpreted in the distributional sense. Here, $X_\sharp\mu$ is the push-forward of measure $\mu$ by diffeomorphism $X$.

In the sequel, we will fix our notation that $A$ will always denote the spatial inverse of diffeomorphism $X$. The following proposition calculates the gradient, Hessian and Gram operator of constraint functions $f_a$ in this setting.
\begin{prop}\label{Geometric Computations of SDiff}
 The following statements hold:
\begin{itemize}
    \item (I) For any $a\in D$, the gradient $\grad f_a$ of $f_a$ is:
    \begin{equation}
     \grad f_a(X)=-\rho_0^{-1}(\nabla^*X)^{-1}\nabla\delta_a=\mathcal N_X(\delta_{X(a)})=\mathcal N_X(\delta_a\circ A)   
    \end{equation}
    \item (II) For any $a\in D$, the Hessian $\Hess f_a$ of $f_a$ is 
    \begin{equation}
    \Hess f_a(X)(\xi,\eta)=-\Tr[\nabla(\xi\circ A)\nabla(\eta\circ A)](X(a)),\quad\forall\xi,\eta\in T_X\MPD    
    \end{equation}
    \item (III) For any $a,a^\prime\in D$, the Gram operator $(Jf)_{a,a^\prime}$ is
    \begin{equation}
     (Jf)_{a,a^\prime}:=\sG_X( \grad f_a(X),\grad f_{a^\prime}(X))=-\Delta_\rho (\delta_{a^\prime}\circ A)(X(a)) .
    \end{equation}
    While the inverse $(Jf)_{a,a^\prime}^{-1}$ reads
    \begin{equation}
      (Jf)_{a,a^\prime}^{-1}(X)=\mathcal G_\rho(X(a),X(a^\prime)) . 
    \end{equation}
    where $\Delta_\rho=\nabla\cdot(\rho^{-1}\nabla)$ and $\mathcal G_\rho$ is the Green's function of the elliptic operator $\Delta_\rho$, under mean-zero condition.
\end{itemize}    
\end{prop}

\begin{proof}[Proof of Proposition~\ref{Geometric Computations of SDiff}]
    The proof follows from straightforward computation. For (I) we have
\begin{equation*}
\sG_X(\grad f_a(X),\delta X)=\frac{d}{d\varepsilon}\leps f_a(X_\varepsilon)=\langle \nabla\cdot(\delta X\circ A)\circ X,\delta_a\rangle_* .   
\end{equation*}
\begin{equation*}
=  -\langle \delta X\circ A,\nabla(\delta_a\circ A)\rangle_*=-\langle\delta X,(\nabla^*X)^{-1}\nabla\delta_a\rangle_*  
\end{equation*}
So we have for any $\delta X\in T_X\MPD$:
\begin{equation*}
  \int_D (\grad f_a(X),\delta X)\rho_0(x)dx= -\int_D \delta X\cdot (\nabla^*X)^{-1}\nabla\delta_a dx.
\end{equation*}
Hence, we conclude
\begin{equation}
 \grad f_a(X)= -\rho_0^{-1} (\nabla^*X)^{-1}\nabla\delta_a 
\end{equation}
which completes the proof of (I). For (II), picking any vector fields $u,v\in \mathfrak X_\sigma$ and denoting $X_\varepsilon:=X+\varepsilon v\circ X$, we have:
\begin{equation*}
\Hess f_a(X)(u\circ X,v\circ X)=\frac{d}{d\varepsilon}\leps \sG_X(\grad f_a(X_\varepsilon),u\circ X).   
\end{equation*}
\begin{equation*}
 = \int_D -\rho_0^{-1}\frac{d}{d\varepsilon}\leps  (\nabla^*X_\varepsilon)^{-1}\nabla\delta_a \cdot u(X)\rho_0(x)dx 
\end{equation*}
\begin{equation*}
 =\int_D (\nabla^*X)^{-1}\nabla^*(v\circ X)(\nabla^*X)^{-1}\nabla\delta_a \cdot u(X)dx .   
\end{equation*}
\begin{equation*}
=\int_D   \nabla\delta_a  (\nabla Au\cdot\nabla v)(X)dx
\end{equation*}
\begin{equation*}
 = -\int_D \delta_a\nabla\cdot [A_\sharp(u\cdot\nabla v)]dx
\end{equation*}
\begin{equation*}
=-\nabla\cdot(u\cdot\nabla v)X(a)=-\mathsf {Tr}(\nabla u\cdot \nabla v)(X(a)).
\end{equation*}
which yields (II). Now for (III), we compute:
\begin{equation*}
\sG_X( \grad f_a(X),\grad f_{a^\prime}(X))=\int_D (\nabla^*X)^{-1}\nabla\delta_a(x)(\nabla^*X)^{-1}\nabla\delta_{a^\prime}(x) \rho_0^{-1}(x)dx   
\end{equation*}
\begin{equation*}
=\int_D   \nabla(\delta_a\circ A)(X)\nabla(\delta_{a^\prime}\circ A)(X)\rho_0^{-1}(x)dx  
\end{equation*}
\begin{equation*}
=\int_D   \nabla(\delta_a\circ A)\nabla(\delta_{a^\prime}\circ A)\rho_0^{-1}\circ Adx 
\end{equation*}
\begin{equation*}
=-\int_D (\delta_a\circ A)\nabla\cdot[\rho^{-1}\nabla(\delta_{a^\prime}\circ A)]    
\end{equation*}
\begin{equation*}
=-\Delta_\rho (\delta_{a^\prime}\circ A)(X(a)) .  
\end{equation*}
Meanwhile, the inverse $(Jf)_{a,a^\prime}^{-1}$ satisfies:
\begin{equation*}
\int_D (Jf)_{a,z}^{-1}(X)(Jf)_{z,a^\prime}(X)dz=\delta(a-a^\prime),    
\end{equation*}
which implies
\begin{equation*}
\delta(a-a^\prime)=-\int_D  (Jf)_{a,z}^{-1}(X)\Delta_\rho (\delta_{a^\prime}\circ A)(X(z))dz   
\end{equation*}
\begin{equation*}
 = -\Delta_\rho [(Jf)_{a,A(z)}^{-1}(X)](X(a^\prime)).
\end{equation*}
Hence, we conclude:
\begin{equation*}
 -\Delta_\rho [(Jf)_{a,A(z)}^{-1}(X)](z)=\delta(X(a)-z) . 
\end{equation*}
As a consequence, we have
\begin{equation*}
 (Jf)_{a,a^\prime}^{-1}(X)=-\mathcal G_\rho(X(a),X(a^\prime)),   
\end{equation*}
where $G_\rho$ is the Green's function of the elliptic operator $\Delta_\rho$.
\end{proof}

\subsection{Geodesic Equation}

We prove in this section the following theorem:
\begin{thm}
 Let $X$ be a geodesic on infinite dimensional Riemannian manifold $(\MPD, \mathsf G)$ with $X(0)=\bId$. Then, the following statements hold:
 \begin{itemize}\label{Geodesic Formalism of IIE}
     \item [(I)] There exists a pressure field $p$ such that
     \begin{equation}\label{IIE Lagrangian ODE}
      \rho_0\ddot X=-\nabla p\circ X.   
     \end{equation}
     \item [(II)] Let $u$ be the velocity field of $X$, then
     \begin{equation}
      p=-\Delta_\rho^{-1}\nabla\cdot(u\cdot\nabla u),\quad \rho=X_\sharp\rho_0  . 
     \end{equation}
     \item[(III)] For $u=\dot X\circ X^{-1}$ and $\rho=X_\sharp\rho_0, $ $(u,\rho,p)$ solves the IIE \eqref{IIE}.
 \end{itemize}
\end{thm}
\begin{proof}[Proof of Theorem~\ref{Geodesic Formalism of IIE}]
    \begin{itemize}
        \item [(I)] By D'Alembert principle, we have $\ddot X\in T_{X}^\perp\MPD$. Hence, there exists a scalar function $\lambda:[0,T]\times  D\to \mathbb R$ such that for any $x\in D$:
\begin{equation}\label{DAlembert}
\rho_0(x)\ddot X(x)=\int_D \lambda(t,z)(\nabla^*X)^{-1}(x)\nabla_z\delta_x(z)dz    
\end{equation}
\begin{equation*}
 =(\nabla^*X)^{-1}(x)\int_D \underbrace{\nabla_z(\delta_z(x)\lambda(t,z))}_{=0\quad a.e.}-\delta_z(x)\nabla\lambda(t,z)dz   
\end{equation*}
\begin{equation*}
 =-(\nabla^*X)^{-1}(x)\int_D   \nabla^*X(z)\nabla(\underbrace{\lambda\circ A}_{:=p})(X(z))\delta_z(x)dz
\end{equation*}
\begin{equation*}
=-(\nabla^*X)^{-1}(x)\int_D \nabla^*X(z)\nabla p(X(z))\delta_x(z)dz=-\nabla p(t,X) ,   
\end{equation*}
\item[(II)] We take a closer look at the Lagrange multiplier $\lambda$. Since for any $a\in D$ we have:
\begin{equation}\label{Level Set Constraint}
 \mathsf G_X(\grad f_a(X),\dot X)=0  . 
\end{equation}
Differentiating \eqref{Level Set Constraint} in time and denoting by $u=\dot X_t\circ X_t^{-1}$, we have for any fixed $a\in D$:
\begin{equation*}
 0=-\frac{d}{dt}\int_D   (\nabla^*X_t)^{-1}\nabla\delta_a(x)\cdot \dot X_t(x)dx  .
\end{equation*}
\begin{equation*}
=-\int_D (\nabla^*X_t)^{-1}\nabla\delta_a(x)\cdot \ddot X_t(x)dx +\int_D (\nabla^*X_t)^{-1}\nabla^*(\dot X_t(x))(\nabla^*X_t)^{-1}\nabla\delta_a\cdot\dot X_t(x)dx   
\end{equation*}
\begin{equation*}
=-\int_D (\nabla^*X_t)^{-1}\nabla\delta_a(x)\cdot \ddot X_t(x)dx +\int_D  \nabla\delta_a(x)(A_t)_\sharp(u\cdot\nabla u)dx  
\end{equation*}
\begin{equation*}
 = -\int_D (\nabla^*X_t)^{-1}\nabla\delta_a(x)\cdot \ddot X_t(x)dx -\int_D  \delta_a(x)[\nabla\cdot(u\cdot\nabla u)](X_t(x))dx  .  
\end{equation*}
Plugging \eqref{DAlembert} into the above equality, we have the following.
\begin{equation*}
 [\nabla\cdot(u\cdot\nabla u)](X_t(a))= -\int_D (\nabla^*X_t)^{-1}(x)\nabla\delta_a(x)  \int_D \rho_0^{-1}(x)\lambda(z)(\nabla^*X_t)^{-1}(x)\nabla_z\delta_x(z)dzdx
\end{equation*}
\begin{equation*}
= -\int_D \int_D \nabla(\delta_a\circ A_t)(X_t(x))\rho_0^{-1}(x)\lambda(z)(\nabla^*X_t)^{-1}(x)\nabla_z\delta_x(z)dzdx 
\end{equation*}
\begin{equation*}
 =\int_D\int_D \nabla(\delta_a\circ A_t) \rho_t^{-1}(x)\nabla^*A_t(x)\nabla\lambda(z)\delta_{A_t(x)}(z)dzdx
\end{equation*}
\begin{equation*}
=\int_D  \nabla(\delta_a\circ A_t) \rho_t^{-1}(x)\nabla(\lambda\circ A_t)(x)dx 
\end{equation*}
\begin{equation*}
 =-  \nabla\cdot (\rho_t^{-1} \nabla p)(X_t(a)),\quad p:=\lambda\circ A_t.
\end{equation*}
Hence, we conclude:
\begin{equation}
  p=-\Delta_\rho^{-1}\nabla\cdot(u\cdot\nabla u),\quad \lambda=p\circ X . 
\end{equation}
\item[(III)] The above geodesic equation \eqref{IIE Lagrangian ODE} in Eulerian coordinates reads:
\begin{equation*}
\rho(\partial_t u+u\cdot\nabla u)-\nabla\Delta_\rho^{-1}\nabla\cdot(u\cdot\nabla u)=0 
\end{equation*}
Together with $\nabla\cdot u=0$ and $\rho$ transported by $u$, we readily verify that IIE \eqref{IIE} is the geodesic equation on $\MPD$ with our weighted metric $\sG$.
\end{itemize}

\end{proof}

\subsection{Second Fundamental Form and Curvature Tensor}
Viewing $\MPD$ as a submanifold of $\Diff$ of infinite dimension and codimension, it makes sense to talk about the second fundamental form of $\MPD$, which admits a concrete description as the projection of the convective acceleration onto the normal direction. To see this, we define for any given $\rho\in\Dens$ the following operator
\begin{equation}
\mathbf{Q}_\rho: w\to \rho^{-1}\nabla\Delta_\rho^{-1}\nabla\cdot w    
\end{equation}
The next lemma shows that $\mathbf Q_\rho$ is an orthogonal projector in weighted $L^2$-sense.

\begin{prop}\label{weighted projector}
 Let $\rho\in\Dens$ and smooth vector field $w\in \mathfrak X(D)$, then the following statements hold:
\begin{itemize}
    \item [(I)] There exists a unique div-free vector field $v$ and $p\in C^\infty(D)$ such that
    \begin{equation}
     w=\rho v+\nabla p.  
    \end{equation}
 moreover, $v=(\bId-\mathbf{Q}_\rho)(\rho^{-1}w)=:\mathbf{P}_\rho(\rho^{-1}w)$. Here, $v$ is not necessarily mean-zero, but $\rho v\in\mathfrak X_0(D)$  iff $w\in \mathfrak X_0(D)$.

 \item [(II)] Let $\mathcal H_\rho:=\{\rho^{-1}\nabla\varphi:\varphi\in C^\infty(D)\}$, then we have $\Ker(\bId-\mathbf{Q}_\rho)=\mathcal H_\rho$. In particular, we have $\mathcal H_\rho\perp \Ran(\bP_\rho)$ in the sense of $L_\rho^2$.

\end{itemize}
\end{prop}
\begin{proof}[Proof of Proposition~\ref{weighted projector}]
\begin{itemize}
    \item [(I)] Given vector field $w$, for any given $\rho\in\Dens$ we write
\begin{equation}\label{Weighted decomposition}
w=\rho\underbrace{[\rho^{-1}w-\rho^{-1}\nabla\Delta_\rho^{-1}\nabla\cdot(\rho^{-1}w)]}_{:={\bP}_\rho(\rho^{-1}w).}+\underbrace{\nabla\Delta_\rho^{-1}\nabla\cdot(\rho^{-1}w)}_{=\rho\mathbf{Q}_\rho(\rho^{-1}w)}
\end{equation}
Let $ v=\bP_\rho(\rho^{-1}w)$ and $p=\Delta_\rho^{-1}\nabla\cdot(\rho^{-1}w)$,
we find the decomposition. Moreover, assume $w=\rho v^\prime+\nabla p^\prime$ for some div-free $v^\prime$ and $p^\prime\in C^\infty(D)$, then $p^\prime$ can be recovered from solving the following elliptic problem:
\begin{equation*}
  \Delta_\rho p^\prime=\nabla\cdot(\rho^{-1}\nabla p^\prime)=\nabla\cdot(\rho^{-1}w)-\underbrace{\nabla\cdot v^\prime}_{=0}=\nabla\cdot(\rho^{-1}w)  .
\end{equation*}
Uniform ellipticity of $\Delta_\rho$ implies the uniqueness of solution $p^\prime$, hence the uniqueness of $v^\prime$.

\item[(II)] Picking any $\varphi\in C^\infty(D)$, we compute $\bP_\rho(\rho^{-1}\nabla\varphi)$ as follows:
\begin{equation*}
\bP_\rho(\rho^{-1}\nabla\varphi)=\rho^{-1}\nabla\varphi- \rho^{-1}\nabla\Delta_\rho^{-1}\nabla\cdot(\rho^{-1}\nabla\varphi)=0   ,
\end{equation*}
which guaranties $\mathcal H_\rho\subset \Ker(\bP_\rho)$. Conversely, assume we have $\bP_\rho \xi=0$, then by the decomposition:
\begin{equation*}
\xi=\mathbf Q_\rho \xi=\rho^{-1}\nabla\Delta_\rho^{-1}\nabla\cdot \xi\in \mathcal H_\rho.
\end{equation*}
which completes the proof of (II).

\end{itemize}
\end{proof} 

The above proposition gives a dynamical decomposition of vector fields into vertical ($\mathfrak X_\sigma$) and horizontal ($\mathcal H_\rho$) directions, from the ambient $\Diff$ viewpoint.

We now connect the second fundamental form $\Pi$ of $\MPD$ w.r.t metric $\mathsf G$ and the above weighted projector $\mathbf Q_\rho$:
\begin{thm}\label{second fundamental form on MPD}
    
The second fundamental form $\Pi_X$ of $\MPD$ at $X$ is given by
\begin{equation}\label{Second Fundamental Form}
\Pi_X(u\circ X,v\circ X)=\mathbf Q_\rho(u\cdot\nabla v),\quad \forall u,v\in \mathfrak X_\sigma(D)\text{ and }X\in\MPD .   
\end{equation}
\end{thm}
\begin{proof}[Proof of Theorem~\ref{second fundamental form on MPD}]
Formula \eqref{Second Fundamental Form} follows from straightforward computation. Mimicking the definition for second fundamental form in finite dimension:
\begin{equation*}
\Pi_X(u\circ X,v\circ X)(x)=\int_{D\times D}(Jf)_{y,z}^{-1}(X)\Hess f_z(X)(u\circ X,v\circ X)\grad f_y(X(x))dydz    
\end{equation*}
\begin{equation*}
=-\int_D[\Delta_\rho^{-1}\nabla\cdot(u\cdot\nabla v)](X(y)) \rho_0^{-1}(x)(\nabla^*X(x))^{-1}\nabla_y\delta_x(y)dy
\end{equation*}
\begin{equation*}
 =-\rho_0^{-1}(x)(\nabla^*X)^{-1}(x)\int_D  [\Delta_\rho^{-1}\nabla\cdot(u\cdot\nabla v)](X(y)) \nabla_y\delta_x(y) dy
\end{equation*}
\begin{equation*}
 =\rho_0^{-1}(x)(\nabla^*X)^{-1}(x)\int_D  \nabla^*X(y)\nabla[\Delta_\rho^{-1}\nabla\cdot(u\cdot\nabla v)](X(y)) \delta_x(y) dy
\end{equation*}
\begin{equation*}
 =\rho_0^{-1}(x)  \nabla[\Delta_\rho^{-1}\nabla\cdot(u\cdot\nabla v)](X(x)) 
\end{equation*}
\begin{equation*}
 =(\rho^{-1}\nabla[\Delta_\rho^{-1}\nabla\cdot(u\cdot\nabla v)])(X(x))= \mathbf Q_\rho (u\cdot\nabla v)(X) .
\end{equation*}
\end{proof}

\begin{rmk}
 $\Pi_X$ is symmetric on $T_X\MPD\times T_X\MPD$, since
\begin{equation*}
\Pi_X(u\circ X,v\circ X)-\Pi_X(v\circ X,u\circ X)=\mathbf{Q}_\rho([u,v])(X)=0    .\end{equation*}
Here, we use the fact that $\nabla\cdot[u,v]=0$ provided $u,v\in\mathfrak X_\sigma(D)$. 
\end{rmk}
\quad

With the above explicit formula for $\Pi$, we can now compute the curvature tensor as follows: For tangent vector $\xi=u\circ X, \eta=v\circ X\in T_X\MPD$ with $u,v\in \mathfrak X_\sigma(D)$, we have
\begin{equation*}
 \mathcal C_{\xi,\eta}=\int_D [\Pi_X(\xi,\xi)(z)\cdot\Pi_X(\eta,\eta)(z)-\lvert \Pi_X(\xi,\eta)(z) \rvert^2] \rho_0(z)dz  
\end{equation*}
\begin{equation*}
 =  \int_D  [\mathbf Q_\rho(u\cdot\nabla u)\mathbf{Q}_\rho(v\cdot \nabla v)-\lvert \mathbf Q_\rho(u\cdot\nabla v)\rvert^2]\rho(z)dz.
\end{equation*}
In view of the above curvature formula,  we recover the Lagrangian instability of shear flows, as steady states of the homogeneous Euler equation.

\section{Hamiltonian and Vorticity Formulation}

In this section, we flip to the Hamiltonian side of the above geometric picture of IIE and derive the Lagrangian-Eulerian formulation of IIE from Hamilton's action principle. The Hamiltonian formulation of the motion of ideal fluids reads coadjoint motion of one-form on $\mathfrak X_\sigma^*(D)$, lies in the framework of general Lie-Poisson equations \cite{MR99}, while the Leray projector essentially becomes the musical isomorphism from cotangent to tangent vectors (see \cite{MW83}\cite{AK98}). Identifying $\mathfrak X_\sigma(D)$ with the space of scalar(vector, resp.) potentials by $\nabla^\perp(\nabla\times$,resp.) in $2(3$-d,resp.), we obtain the vortex dynamics. 

Here, we extend the above geometric picture to the inhomogeneous setting. A Hamilton-Pontryagin type variational principle and a companion Lagrangian-Eulerian formulation is established, which generalizes the Weber's formula for perfect fluids in \cite{Co01a}. We then establish the vorticity formulation for IIE in simply connected domain and torus respectively. The reason why we need to treat torus case separately is because the mean value of $\rho u$ is conserved, instead of $u$, hence a global stream function representation for $u$ is not applicable in torus case. A proof of local existence of strong solutions for IIE is obtained based on our vorticity formulation.

\subsection{Hamilton-Pontryagin Principle and Lagrangian-Eulerian Formulation}

We now lift ourselves to the cotangent bundle of $\MPD$ where the Hamiltonian equation naturally lives and we will derive it in our setting. To this end, we first write down the kinetic energy as a functional on tangent bundle:
\begin{equation}\label{Lagrangian Action of IIE}
L: (X,\dot X)\in T_X\MPD\to L(X,\dot X):=\frac{1}{2}\int_D \lvert \dot X(t,x)\rvert^2 \rho_0(x)dx     
\end{equation}
It is not hard to show that the action-minimizing problem of $L$, augmented by the incompressibility constraint $\det(\nabla X)=1$, is precisely solved by the solution of equation \eqref{IIE}. We refer the interested readers to \cite{[LLP11]} for study of a relaxed geodesic problem in the sense of Brenier's generalized flow (cf. \cite{Br89}, \cite{Br99}).

Abusing the notation a little bit, we now write the Lagrangian $L$ in an equivalent form:
\begin{equation}
L(u,X)=\frac{1}{2}\int_D \lvert u\rvert^2 (X_\sharp\rho_0)(x)dx     
\end{equation}
Let $D=\mathbb R^d$ or $\mathbb T^d$, we now state the following:
\begin{thm}[(Hamilton-Pontryagin Principle)]\label{H-P Principle}
Fix $T>0$ and $\rho_0\in\Dens$ satisfying assumption \ref{Assumption A}, we assume the triplet 
\begin{equation*}
(X,u,\xi)\in C^1([0,T];\MPD\times\mathfrak X_\sigma(D)\times\mathfrak X(D))    \end{equation*} 
is a critical point of the augmented action
\begin{equation}
\Phi(u,X,\xi)=\int_0^T L(u,X)-\langle\xi,\dot X-u(X)\rangle_{L_{\rho_0}^2} dt .   
\end{equation}
Then the following statements hold:
\begin{itemize}
    \item [(I)]$(u_t,\rho_0\circ X_t^{-1})$ solves the incompressible inhomogeneous Euler equation \eqref{IIE}.

    \item [(II)] Consider the following controlled Hamiltonian:
    \begin{equation*}
     H(\xi,X,u)=L(u,X)+\langle\xi,u(X)\rangle_{L_{\rho_0}^2} .  
    \end{equation*}
    Then $(X,u,\xi)$ solves the associated Hamiltonian ODE:
    \begin{equation}\label{Hamiltonian ODE for IIE}
    \dot X=\frac{\partial H}{\partial\xi},\quad \dot \xi=-\frac{\partial H}{\partial X},\quad \frac{\partial H}{\partial u}=0  .  
    \end{equation}
    \item [(III)] Given the initial data $(u_0,\xi_0)$ satisfying compatibility condition
    \begin{equation*}
        u_0=\bP_{\rho_0}\xi_0,\quad\text{equivalently }\xi_0= u_0+\rho_0^{-1}\nabla q_0.
    \end{equation*}
    the solution of \eqref{Hamiltonian ODE for IIE} is given by the following Lagrangian-Eulerian system:
\begin{equation}\label{IIE Lagrangian}
\left\{
\begin{aligned}
&\dot X_t=u(t,X_t), A_t=X_t^{-1};\\
&\rho(t,\cdot)=\rho_0\circ A_t;\\
& u_t=\bP_\rho\bigg[\nabla^*A_t\bigg(\xi_0-\frac{\nabla \log\rho_0}{2}\int_0^t  \lvert \dot X_\tau(x)\rvert^2d\tau\bigg)\circ A_t\bigg] .
\end{aligned}
\right.
\end{equation}
\end{itemize}
\end{thm} 
\begin{proof}[Proof of Theorem~\ref{H-P Principle}]
\begin{itemize}
    \item [(I)] As a reminder of our notation, we denote by $A=X^{-1}$ and $\rho:=\rho_0\circ A$. We compute the first variation of $\Phi$. We consider smooth perturbations $\delta u\in C^1([0,T];\mathfrak X_\sigma(D))$, $\delta X\in C_c^1([0,T];\MPD)$ and $\delta\xi\in C^1([0,T];\mathfrak X(D))$ of $\Phi$, respectively, and compute the corresponding variational equations as follows:
\begin{itemize}
\item [(i)] The variation in the variable $u$ is given by
\begin{equation}\label{variation of u}
\frac{d}{d\varepsilon}\bigg\rvert_{\varepsilon=0}\Phi(u+\varepsilon\delta u,X,\xi)=\int_0^T \int_{D} \langle u\rho,\delta u\rangle -\langle \xi, \delta u(X)\rangle_{L_{\rho_0}^2} dxdt    .
\end{equation}
Hence, vanishing of \eqref{variation of u} implies that $\rho(u-\xi\circ A)$ is a gradient, i.e. there exists some scalar function $q$ such that
\begin{equation}\label{Stationary u}
 u\rho=(\rho_0\xi)\circ A+\nabla q  .
\end{equation}

\item [(ii)] The variation in variable $X$ is given by
\begin{equation*}
\frac{d}{d\varepsilon}\bigg\rvert_{\varepsilon=0}\Phi(u,X+\varepsilon\delta X,\xi)    
\end{equation*}
\begin{equation}\label{variation in X}
 =\int_0^T \int_{D} \frac{d}{d\varepsilon}\bigg\rvert_{\varepsilon=0}\frac{1}{2}\lvert u\rvert^2\rho_0\circ X_\varepsilon^{-1}+\langle \xi,\dot X+\varepsilon\dot{\delta X}\rangle_{L_{\rho_0}^2}  -\langle \xi, u(X+\varepsilon \delta X)\rangle_{L_{\rho_0}^2} dxdt .   
\end{equation}
To compute $\frac{d}{d\varepsilon}X_\varepsilon^{-1}$, notice:
\begin{equation*}
 0=\frac{d}{d\varepsilon}(X_\varepsilon\circ X_\varepsilon^{-1})=\delta X\circ X_\varepsilon^{-1}+\nabla X_\varepsilon\circ X_\varepsilon^{-1}\cdot \frac{d}{d\varepsilon}X_\varepsilon^{-1}   .
\end{equation*}

Hence, we have:
\begin{equation}\label{Variation in X^-1}
  \frac{d}{d\varepsilon}\bigg\rvert_{\varepsilon=0} X_\varepsilon^{-1}=-[(\nabla X_\varepsilon)^{-1}\delta X]\circ X_\varepsilon^{-1} \bigg\rvert_{\varepsilon=0}=-[(\nabla X)^{-1}\delta X]\circ A.  
\end{equation}
Therefore, combining \eqref{variation in X} and \eqref{Variation in X^-1} we conclude
\begin{equation*}
    \frac{d}{d\varepsilon}\bigg\rvert_{\varepsilon=0}\Phi(u,X+\varepsilon\delta X,\xi)
    \end{equation*}
   \begin{equation*} 
    =\int_0^T \int_{D} -\frac{1}{2}\lvert u\rvert^2\nabla\rho_0\circ A\cdot [(\nabla X)^{-1}\delta X]\circ A+\langle \xi,\dot{\delta X}\rangle_{L_{\rho_0}^2} -\langle \xi, \nabla u(X)\delta X\rangle_{L_{\rho_0}^2} dxdt  
\end{equation*}
\begin{equation}\label{Stationary condition in X}
=-\int_0^T \int_{D} \langle \rho_0\dot \xi+\rho_0\nabla^*u(X)\xi+\frac{1}{2}\lvert u\circ X\rvert^2\nabla\rho\circ X,\delta X\rangle  dxdt    .
\end{equation}
By arbitrariness of $\delta X$ in \eqref{Stationary condition in X}, the variational equation in $X$ reads
\begin{equation}\label{Hamiltonian ODE}
\rho_0[\dot \xi+\nabla^*u(X)\xi]+\frac{1}{2}\lvert u\circ X\rvert^2\nabla\rho\circ X=0  .
\end{equation}
\item [(iii)] The variation in variable $\xi$ is simply given by
\begin{equation}\label{Stationary xi}
\dot X=u(X) .
\end{equation}
Hence, the critical point equation is obtained by putting together \eqref{Stationary u}\eqref{Stationary condition in X}\eqref{Stationary xi}:
\begin{equation}\label{Critical Point Equations}
\left\{
\begin{aligned}
&\dot X=u(X),\quad A=X^{-1};\\
&\rho_0[\dot \xi+\nabla^*u(X)\xi]+\frac{1}{2}\lvert u\circ X\rvert^2\nabla\rho\circ X=0;  \\
& u\rho=(\rho_0\xi)\circ A+\nabla q,\quad \nabla\cdot u=0;\\
&\rho=X_\sharp\rho_0=\rho_0\circ A.
\end{aligned}
\right.
\end{equation}
\end{itemize}
Now denote by $\eta=(\rho_0\xi)\circ A$, by method of characteristics, \eqref{Hamiltonian ODE} implies
\begin{equation*}
\partial_t\eta=(\rho_0\dot\xi)\circ A+\nabla(\rho_0\xi)\circ A\cdot \partial_tA=-\nabla^*u(\rho_0\xi)\circ A-\frac{1}{2}\lvert u\rvert^2\nabla\rho-u\cdot\nabla[(\rho_0\xi)\circ A]   . 
\end{equation*}
hence $\eta$ solves the PDE:
\begin{equation}\label{Dual Lie Transport Eta}
  \partial_t\eta+u\cdot\nabla\eta+\nabla^*u\cdot\eta+\frac{1}{2}\lvert u\rvert^2\nabla\rho=0  .
\end{equation}
Now, since $u\rho=\nabla q+\eta$, plugging in the equation \eqref{Dual Lie Transport Eta} we obtain
\begin{equation*}
\partial_t(u\rho-\nabla q)+u\cdot \nabla(u\rho-\nabla q)+\nabla^*u\cdot(u\rho-\nabla q)+\frac{1}{2}\lvert u\rvert^2\nabla\rho=0  .  
\end{equation*}
Simple rearrangement yields
\begin{equation*}
\partial_t(\rho u)+u\cdot \nabla (\rho u)+\nabla\bigg(\frac{1}{2}\lvert u\rvert^2\rho-u\cdot\nabla q-\partial_t q\bigg)=0  .  
\end{equation*}
Hence, we let
\begin{equation*}
p=\frac{1}{2}\lvert u\rvert^2\rho-u\cdot\nabla q-\partial_t q    
\end{equation*}
and conclude that \eqref{Critical Point Equations} is equivalent to the system \eqref{IIE} with pressure $p$.

\item[(II)] We compute the partial derivatives $\partial_\xi H$, $\partial_XH$. Recall the definition of variational derivative:
\begin{equation*}
 \frac{d}{d\varepsilon}\leps H(\xi+\varepsilon\delta\xi,X,u)= \int_D [\frac{\partial H}{\partial\xi}(\xi,X,u)\cdot\delta\xi]\rho_0(a)da , 
\end{equation*}
and
\begin{equation*}
 \frac{d}{d\varepsilon}\leps H(\xi,X+\varepsilon\delta X,u)=\sG_X\big(\frac{\partial H}{\partial X},\delta X\big)= \int_D [\frac{\partial H}{\partial X}(\xi,X,u)\cdot\delta X]\rho_0(a)da; 
\end{equation*}
which implies the following Hamilton's ODE
\begin{equation*}
 \frac{\partial H}{\partial\xi}=u(X),\quad  \frac{\partial H}{\partial X}=\nabla^*u(X)\xi+\frac{1}{2}\lvert u(X)\rvert^2 (\nabla^*X)^{-1}\nabla\log\rho_0  ;
\end{equation*}
with control
\begin{equation*}
\frac{\partial H}{\partial u}=u\rho-(\rho_0\xi)\circ A-\nabla q=0 .
\end{equation*}
Putting above together, straightforward computation yields that \eqref{Hamiltonian ODE for IIE} is equivalent to \eqref{Critical Point Equations}.

\item[(III)] We consider the following general version of \eqref{Hamiltonian ODE}:
\begin{equation*}
 \dot w+\nabla^*u(X)w=f(X),\quad w(0,\cdot)=w_0 .  
\end{equation*}
Provided sufficient regularity on $X$, applying Duhamel's formula we have the following Lagrangian formula for solution $w(t,\cdot)$:
\begin{equation*}
w(t,x)= (\nabla^*X_t)^{-1}\bigg( w_0(x)+\int_0^t (\nabla^*X_\tau) f(\tau,X_\tau(x))d\tau\bigg)   
\end{equation*}
Hence, we obtain the following representation for the  solution $\eta$ of \eqref{Dual Lie Transport Eta}:
\begin{equation*}
 \eta(t,x)=  \nabla^*A_t\bigg(\rho_0\xi_0-\frac{1}{2}\int_0^t (\nabla^*X_\tau) \lvert u(\tau,X)\rvert^2\nabla \rho(\tau,X_\tau(x))d\tau\bigg)\circ A_t(x)
\end{equation*}
\begin{equation}\label{Lagrangian Formulae of eta}
 = \nabla^*A_t\bigg(\rho_0\xi_0-\frac{1}{2}\int_0^t  \lvert \dot X_\tau(x)\rvert^2\nabla \rho_0(x)d\tau\bigg)\circ A_t(x).  
\end{equation}
 We are now ready to conclude:
\begin{equation}\label{Lagrangian Formulae of u}
u=\bP_\rho(\rho^{-1}\eta)=\bP_\rho\bigg[\nabla^*A_t\bigg(\xi_0-\frac{\nabla \log\rho_0}{2}\int_0^t  \lvert \dot X_\tau(x)\rvert^2d\tau\bigg)\circ A_t\bigg]   .
\end{equation}
Combining \eqref{Stationary xi} \eqref{Lagrangian Formulae of u} and $\rho=\rho_0\circ A$, we recover the system \eqref{IIE Lagrangian}. 
\end{itemize}
\end{proof}

\subsection{Vorticity Formulation}

 We now move on to derive a stream-vorticity formulation of \eqref{IIE}. We will state a more general functional analytic framework and treat IIE as a specific example. Here we treat the case that $D$ is simply connected and $D=\mathbb T^2$ separately.

Let $\mathcal V$ be a Hilbert space with inner product $\langle\cdot,\cdot\rangle_{\mathcal V}$. Let $\mathcal I:\mathcal V\to\mathcal V^*$ be the Riesz isomorphism between dual pairing $\langle \cdot,\cdot\rangle_*$ of $\mathcal V,\mathcal V^*$:
\begin{equation*}
 \langle v,v^\prime\rangle_*=\langle v,\mathcal I^{-1}v^\prime\rangle_{\mathcal V}   
\end{equation*}
 Let $\mathcal D: \mathcal F\to \mathcal V$ be a linear operator from Hilbert space $\mathcal F$ (with dual pairing $\langle,\rangle$) to the Hilbert space $\mathcal V$ with closed range, then we have the following pull-back inner product on $\mathcal F$:
\begin{equation}
 \langle \varphi,\psi\rangle_{\mathcal F}:=\langle \mathcal D\varphi,\mathcal D\psi\rangle_{\mathcal V}=\langle \mathcal D\varphi,\mathcal I\mathcal D\psi\rangle_*=\langle\varphi,\mathcal D^*\mathcal I\mathcal D\psi\rangle
\end{equation}
Denote by $\mathcal A:\mathcal F\to\mathcal F^*$ the Riesz isomorphism on $(\mathcal F,\langle,\rangle_{\mathcal F})$ we have for $\psi^\prime\in\mathcal F^*$ with $\mathcal A\psi=\psi^\prime$:
\begin{equation}
\langle \varphi,\mathcal A\psi\rangle=\langle\varphi,\mathcal D^*\mathcal I\mathcal D\psi\rangle  . 
\end{equation}
Hence, we conclude: $\mathcal A=\mathcal D^*\mathcal I\mathcal D$. While the projection onto the range $\mathsf{Ran}(\mathcal D)$ is:
\begin{equation*}
\mathcal P=\mathcal D\mathcal A^{-1}\mathcal D^*\mathcal I : \mathcal V\to \mathsf{Ran}(\mathcal D)   
\end{equation*}
In our hydrodynamic setting, assume $d=2$ and $D\subset\mathbb R^2$ is a simply connected smooth domain. We let $\mathcal V=L^2(D;\mathbb R^d)$ be the space of $L^2$ vector fields, $\mathcal F=H^{1}(D)/\mathbb R$, $\mathcal F^*=H_0^{-1}(D)$ be the space of $H^{-1}$ distributions. Set $\mathcal D=\nabla^\perp$, hence $\mathcal D^*=\nabla\times$. The inertia operator in this case is the following uniformly elliptic operator (due to assumption \ref{Assumption A}):
\begin{equation*}
 \mathcal A=\nabla\times(\rho\nabla^\perp)=:\bL_\rho=\nabla\cdot(\rho\nabla) ,  
\end{equation*}
and the projector reads:
\begin{equation*}
\mathcal P= \nabla^\perp\bL_\rho^{-1}\nabla\times(\rho (\cdot) ).   
\end{equation*}
Here, the inverse $\bL_\rho^{-1}$ is understood in the sense of Dirichlet boundary condition.

Indeed, decomposing $\rho^{-1}w=v+\rho^{-1}\nabla p$ as in proposition 2.5(I), we have the following identity
\begin{equation}\label{Weighted Biot Savart}
 \mathcal P(\rho^{-1}w)=\nabla^\perp\bL_\rho^{-1}\nabla\times(\rho v)=:\mathcal K_\rho \nabla\times(\rho v)   .
\end{equation}
Here, $\mathcal K_\rho$ is the modified Biot-Savart operator, which plays the role of the right inverse of the $\rho$-weighted curl operator.  \eqref{Weighted Biot Savart} is a generalization of the following identity in the homogeneous setting:  

\begin{equation*}
 \bP v=\nabla^\perp\Delta^{-1}\nabla\times v,\quad \forall v\in\mathfrak X(D).   
\end{equation*}

Since $d=2$ and $\nabla\cdot v=0$, let $\psi$ be the stream function of $v$, we find:
\begin{equation*}
 \mathcal P(\rho^{-1}w)=\nabla^\perp\bL_\rho^{-1}\nabla\times(\rho \nabla^\perp\psi)= \nabla^\perp\bL_\rho^{-1}\bL_\rho \psi=\nabla^\perp\psi=v.  
\end{equation*}
Therefore $\mathcal P=\bP_\rho$. Similar calculation shows that $\mathcal P=\bP_\rho$ also holds in $d=3$, if we replace $\nabla^\perp$ with $\nabla\times$ and define $\bL_\rho$ as the following weighted Laplacian for vector fields:
\begin{equation*}
\bL_\rho: w\in\mathfrak X(D)\to \nabla\times(\rho\nabla\times w).    
\end{equation*}

We \textit{formally} unify the above relations in  $2$-d in the following commutative diagram.

\begin{equation}
\begin{tikzcd}[row sep=large, column sep=large]
{\mathcal F := H_0^{1}(D)}
  \arrow[r, "{\mathcal D=\nabla^\perp}"]
&
{\mathcal V := L^2_{\rho}(D;\mathbb{R}^{2})}
  \arrow[d, "{\mathcal I:\ v\mapsto \rho v}"]
  \arrow[loop right, distance=2.2em, out=20, in=-20, looseness=6,
                  "{\mathbf{P}_{\rho}=\mathcal{D}\mathcal{A}^{-1}\mathcal{D}^*\mathcal{I}}"]
\\
{\mathcal F^{*} := H_0^{-1}(D)}
\arrow[u,"{\mathcal A^{-1}=\mathbf{L}_\rho^{-1}}"]
\arrow[ur,"\mathcal K_\rho=\nabla^\perp\bL_\rho^{-1}"] 
&
{\mathcal V^{*} \cong L^2_{\rho}(D;\mathbb{R}^{2})}
\arrow[l,"{\mathcal D^*=\nabla\times}"]
\end{tikzcd}
\end{equation}

Following the geometric picture above, we now introduce the vorticity formulation in this simply-connected setting. Let $\eta=\rho u+\nabla q$, we therefore have $\bP_\rho (\rho^{-1}\eta)=u$. Following our above computation, we have yet another closed-form representation of $\bP_\rho$:
\begin{equation*}
\bP_\rho(\rho^{-1}\eta)=\nabla^\perp\bL_\rho^{-1}\nabla\times\eta .   
\end{equation*}

Now we define $\omega=\nabla\times \eta$, then taking the curl of \eqref{Hamiltonian ODE for IIE} we obtain the following vorticity equation:
\begin{equation}
 \partial_t\omega+\mathscr L_u \omega+\nabla\frac{1}{2}\lvert u\rvert^2\times\nabla \rho=0  . 
\end{equation}
where $\mathscr L_u$ is the Lie derivative of a scalar/vector field in $2d/3d$ respectively. Hence, in view of the above diagram, we have the following vorticity formulation of \eqref{IIE}, where the pressure is deleted:
\begin{equation}\label{Vorticity IIE}
\left\{
\begin{aligned}
&\partial_t\omega+\mathscr L_u \omega+\nabla\frac{1}{2}\lvert u\rvert^2\times\nabla \rho=0\\
&u=\mathcal K_\rho \omega:=\nabla^\perp\bL_\rho^{-1}\omega\\
&\partial_t\rho+u\cdot\nabla\rho=0
\end{aligned}
\right.
\end{equation}
We note that \eqref{Vorticity IIE} is algebraically equivalent to system \eqref{IIE Lagrangian}. The boundary condition is $u\cdot \mathbf{n}=0$ on $\partial D$, which is compatible with the Dirichlet boundary condition of the Poisson equation 
\begin{equation}\label{Elliptic PDE for Stream}
 \bL_\rho \psi=\omega,\quad \psi \big\rvert_{\partial D}=0.   
\end{equation}

We now show the conservation of vorticity (Kelvin circulation theorem) in the framework of IIE. In hydrodynamics, conservation laws coming from particle relabeling symmetry (Casimir Invariants in geometry) are known to be the consequence of duality between Lie transport (Lie elements of Lagrangian loops) and dual Lie transport (co-adjoint motion of momentum, cf.\cite{HMR98}\cite{MP25}). Here, for the inhomogeneous fluid, we do not have a full conservation of circulation, due to the forcing term in vortex dynamics. However, the following result holds for the case that $D\subset \mathbb R^2$ is simply connected:
\begin{thm}[(Conservation of Circulation)]\label{Kelvin Circulation IIE}
Assume $(u,\rho)$ is a smooth solution of IIE. For any regular value $\alpha> 0$ of $\rho_0$, assume $E_\alpha\subset\{\rho_0=\alpha\}$ is a smooth loop, let $D_\alpha:=\{\rho_0<\alpha\}$ be the smooth bounded domain enclosed by $E_\alpha$, then we have
\begin{equation}\label{Circulation of IIE}
\oint_{X_t(E_\alpha)} u_t\cdot d\ell=\oint_{E_\alpha} u_0\cdot d\ell,   
\end{equation}
which implies 
\begin{equation}\label{Conserved Vorticity}
 \int_{X_t(D_\alpha)}\omega_t(x) dx=\int_{D_\alpha}\omega_0(x)dx .  
\end{equation}
\end{thm} 
\begin{proof}[Proof of Theorem~\ref{Kelvin Circulation IIE}]
Thanks to the Lagrangian formula \eqref{Lagrangian Formulae of eta}, we have for any smooth loop $\gamma$:
\begin{equation*}
\oint_{X_t(\gamma)}\eta_t\cdot d\ell=\oint_\gamma \nabla^*X_t\eta_t\circ X_t \cdot d\ell 
\end{equation*}
\begin{equation*}
=\oint_\gamma \big[u_0\rho_0+\nabla q_0-\nabla\rho_0(x)\frac{1}{2}\int_0^t\lvert \dot X_\tau(x)\rvert^2 d\tau   \big]\cdot d\ell    
\end{equation*}
\begin{equation*}
 =\oint_\gamma u_0\rho_0\cdot d\ell-\frac{1}{2}\int_0^t\oint_\gamma\lvert \dot X_\tau\rvert^2 \nabla\rho_0  \cdot d\ell d\tau .  
\end{equation*}
Now choose $\gamma=E_\alpha$, since $\{\rho_t=\alpha\}=X_t(E_\alpha)$, we conclude:
\begin{equation*}
 \oint_{X_t(E_\alpha)} u_t\cdot d\ell= \alpha^{-1}\oint_{X_t(E_\alpha)} \eta_t\cdot d\ell=\alpha^{-1}\oint_{E_\alpha}\rho_0 u_0\cdot d\ell=\oint_{E_\alpha}u_0\cdot d\ell .
\end{equation*}
where we use the orthogonality relation: For all smooth $f$ we have
\begin{equation*}
\oint_{E_\alpha} f\nabla\rho_0\cdot d\ell=\int f\nabla\rho_0\cdot \frac{\nabla^\perp\rho_0}{\lvert\nabla\rho_0\rvert}\mathscr H_{E_\alpha}^1(dx)=0.    
\end{equation*}
Hence, the proof of \eqref{Circulation of IIE} is complete. To show \eqref{Conserved Vorticity}, simply apply the Stokes theorem to \eqref{Circulation of IIE} multiplied by $\alpha$.

\end{proof}

\begin{rmk}[(Noether Theorem)] 
 Actually, we may derive the conservation of circulation of IIE as a consequence of Noether's theorem. By \eqref{Lagrangian Action of IIE}, we consider the following action:
 \begin{equation*}
  S(\{X_t\})=\int_0^T L(X,\dot X)dt=\int_0^T\int _D \frac{1}{2}\lvert\dot X(t,a)\rvert^2\rho_0(a)da   
 \end{equation*}
 We notice that $S$ is invariant under right action of $\Psi\in \mathsf{Diff}_{\rho_0}(D)$, the subgroup of $\Diff$ preserving the initial density $\rho_0$. Moreover, by Lemma \ref{Abel-Jacobi-Liouville Formula}, for curve $\Psi_\varepsilon\subset \mathsf{Diff}_{\rho_0}(D)$ generated by vector field $v$ and starting from identity, we have:
 \begin{equation*}
 0=\frac{d}{d\varepsilon}\bigg\rvert_{\varepsilon=0}\rho_0=\frac{d}{d\varepsilon}\bigg\rvert_{\varepsilon=0}\rho_0(\Psi_\varepsilon) \det(\nabla\Psi_\varepsilon)=\rho_0\nabla\cdot v +v\cdot\nabla\rho_0=\nabla\cdot(v\rho_0)
 \end{equation*}
 Hence, $\rho_0$-preserving flows are generated by vector field $v$ such that $\nabla\cdot(v\rho_0)=0$. Now, combined with the $\MPD$-constraint of the IIE dynamics, we conclude that $\mathscr G=\MPD\cap \mathsf{Diff}_{\rho_0}(D)$ is the right gauge group for the IIE action $S$, which is significantly smaller than $\MPD$ and fails the full Kelvin circulation theorem. In view of the above computation, we have:
 \begin{equation*}
 T_X \mathscr G=\{v\circ X: \nabla\cdot v=v\cdot\nabla\rho_0=0\}.   
 \end{equation*}
Differentiating $S$ with right action of $\{\Psi_\varepsilon\}\subset \mathscr G$, we compute:
 \begin{equation*}
  \frac{d}{d\varepsilon}\bigg\vert_{\varepsilon=0} S(\{X_t\circ \Psi_\varepsilon\})=  \frac{d}{d\varepsilon}\bigg\vert_{\varepsilon=0}\int_0^T\int _D \frac{1}{2}\lvert\dot X(t,\Psi_\varepsilon(a))\rvert^2\rho_0(a)da    
 \end{equation*}
 \begin{equation*}
  = \int_0^T\int _D \dot X(t,a)[\nabla \dot X(t,a)v(a)]\rho_0(a)dadt
 \end{equation*}
 \begin{equation*}
  =   \int_0^T \int_D \frac{d}{dt} [\dot X(t,a)(\nabla X(t,a)v)]\rho_0(a)dadt-\int_0^T\int_D \rho_0(a)\ddot X(t,a)[\nabla X(t,a)v]da
 \end{equation*}
 \begin{equation*}
  =  \int_D u(T,a)\cdot(X_T)_\sharp(v\rho_0) da-\int_D u_0(a)\cdot (v\rho_0)da-\int_0^T\int_D \nabla p(t,a) (X_\sharp v)(a)da
 \end{equation*}
 \begin{equation*}
 =  \int_D u(T,a)\cdot(X_T)_\sharp(v\rho_0) da-\int_D u_0(a)\cdot (v\rho_0)da+\int_0^T\int_D p(t,a)  (\nabla\cdot v)\circ X^{-1}(a)da    
 \end{equation*}
 \begin{equation*}
  =  \int_D u(T,a)\cdot(X_T)_\sharp(v\rho_0) da-\int_D u_0(a)\cdot (v\rho_0)da   
 \end{equation*}
 where we integrate by part, use \eqref{IIE Lagrangian ODE}, change the variable and apply the following commutation relation between push-forward and divergence:
 \begin{equation*}
  \nabla\cdot(X_\sharp v)=(\nabla\cdot v)\circ X^{-1}.   
 \end{equation*}
 Hence, we conclude for all $v\in \tilde{\mathfrak g}:=T_e\mathscr G$, the quantity
 \begin{equation}\label{Noether Quantity}
  \mathscr I_t:=\int_D u(t,a)\cdot [X_\sharp (\rho_0v)](t,a) da  
 \end{equation}
 is conserved. Typical example includes the case $d=2$ and $v=F(t,\rho_0)\nabla^\perp \rho_0$ for any smooth function $F: \mathbb R_+\times\mathbb R_+\to\mathbb R$. Plugging such $v$ back to \eqref{Noether Quantity}, we have
 \begin{equation*}
  \mathscr I_t=\int_D  u(t,a) \nabla^\perp\rho(t,a) \rho(t,a) F(\rho(t,a))da,\quad \rho(t,a)=\rho_0(X_t^{-1}(a)).
 \end{equation*}
 Assume that $\Ran(\rho_t)=[0,c]$ (here $c$ is independent of $t$ since $\rho_t$ is a rearrangement of $\rho_0$), by coarea formula we have:
 \begin{equation*}
 \mathscr I_t=\int_0^c \alpha F(\alpha)\bigg(\int_{\{\rho_t=\alpha\}}  u(t,a)\cdot\frac{\nabla^\perp \rho_t}{\lvert\nabla \rho_t\rvert} d\mathscr H^1\bigg)d\alpha=\int_0^c \alpha F(\alpha)  \oint_{\{\rho_t=\alpha\}}  u(t,a)\cdot d\ell d\alpha
 \end{equation*}
 \begin{equation*}
 =    \int_0^c \alpha F(\alpha)  \Gamma(t,\alpha) d\alpha,\quad \Gamma(t,\alpha):=\oint_{\{\rho_t=\alpha\}}  u(t,a)\cdot d\ell .
 \end{equation*}
 Thanks to the conservation of $\mathscr I_t$ and arbitrariness of $F$, we conclude $\Gamma(t,\alpha)=\Gamma(0,\alpha)$ for all $t$ and $\alpha\in \Ran(\rho_t)$.
\end{rmk}

\begin{rmk}
 It is a natural question to ask that whether such vorticity formulation would lead the way to the study of finite time blow up of IIE in dimension $2$ and $3$. We remark that the Boussinesq system, which is a fluid equation with active tracers coupled via Buoyancy, is demonstrated in the recent breakthrough \cite{ElgindiPasqualotto2023}\cite{[EP25]} to blow up in finite time in its local well-posedness class. However, in view of the vortex dynamics, the source term in IIE is much harder to control compared to the Boussinesq case. On the other hand, analysis of lifespan and continuation criteria in IIE based on the vorticity formulation is applicable, see the companion paper \cite{Pa26+}.   
\end{rmk}

We now discuss the stream-vorticity formulation in torus. Assume $D=\mathbb T^2$, the key observation is that: For IIE, the conservation of linear momentum reads
\begin{equation*}
 \frac{d}{dt}\int_D   \rho u dx=-\int_D  u\cdot\nabla (\rho u)dx-\int_D \nabla p dx=0.
\end{equation*}
Which follows from integrating the first term by part and using mean-free property of gradient fields on torus. However, mean of $u$ is not conserved.

By the above observation, WLOG we assume that $\rho_0u_0$ is mean zero. Abuse our notation a little bit temporarily, we let
\begin{equation*}
 \eta=\rho u=\nabla^\perp\psi+\nabla q  . 
\end{equation*}
Now, the dynamical equation for $\omega$ and $\rho$ keep unchanged:
\begin{equation*}
 \partial_t\omega+u\cdot\nabla\omega+\{\frac{1}{2}\lvert u\rvert^2,\rho\}=0,\quad \partial_t\rho+u\cdot\nabla\rho=0,   
\end{equation*}
where $\{\cdot,\cdot\}$ is the standard Poisson bracket in $2$-d.

The constitutive law becomes different: Taking curl on $\eta$ and making use of $\nabla\cdot u=0$, we have:
\begin{equation*}
  \omega=\Delta\psi,\quad  u=\rho^{-1}\nabla^\perp\psi+\rho^{-1}\nabla p\quad\Longrightarrow\quad u= \bP_\rho(\rho^{-1}\nabla^\perp\Delta^{-1}\omega).
\end{equation*}
Hence, the stream-vorticity formulation of IIE on $\mathbb T^2$ reads:
\begin{equation}\label{Vorticity IIE on Torus}
\left\{
\begin{aligned}
&\partial_t\omega+\mathscr L_u \omega+\{\frac{1}{2}\lvert u\rvert^2,\rho\}=0;\\
&u=\bP_\rho(\rho^{-1}\nabla^\perp\Delta^{-1}\omega);\\
&\partial_t\rho+u\cdot\nabla\rho=0.
\end{aligned}
\right.
\end{equation}
Supplemented by mean-zero condition of $\rho_0 u_0$. Notice that the above constitutive law is a general law and reduces to the weighted Biot-Savart law in case that $D$ is a simply connected domain:
\begin{equation*}
 \bP_\rho(\rho^{-1}\nabla^\perp\Delta^{-1}\omega) =\mathcal K_\rho \nabla\times(\rho\rho^{-1}\nabla^\perp\Delta^{-1}\omega)=\mathcal K_\rho\omega=\nabla^\perp\bL_\rho^{-1}\omega .  
\end{equation*}

\subsection{Local Existence Via Vorticity Formulation}

We recall that in the general dimension $d$, IIE is shown to be locally well-posed in $H^{s+1}\times H^{s+1}$ for $s>d/2$, and we refer the readers to \cite{Da06}\cite{[Da10]} \cite{[DF11]} for various well-posedness results in $L^p$-based Sobolev spaces and endpoint Besov spaces. 

Here, we show the local existence of strong solutions of IIE in Sobolev space via the vorticity formulation \eqref{Vorticity IIE}. We choose to work on a  smooth and simply connected domain $D\subset \mathbb R^2$ for simplicity, but the proof of local well-posedness on $\mathbb T^2$ using \eqref{Vorticity IIE on Torus} is also available following the same line of argument. We also expect that the Lagrangian–Eulerian system \eqref{IIE Lagrangian} yields a unified proof of local well-posedness for strong solutions in all dimensions. For a similar Lagrangian approach to the local well-posedness of hydrodynamic models, we refer the interested reader to \cite{Co01a}\cite{CI08}\cite{MP25}\cite{Zh10}. However, we will only focus on the proof through our vorticity formulation to highlight its utility.

We start with the following elementary elliptic lemma. The proof is standard and is therefore omitted.
\begin{lem}\label{Ellptic Regularity}
Assume that $\omega\in H^2$ and $\rho\in H^{3}$ satisfies the assumption \ref{Assumption A}, then the solution $\psi$ of \eqref{Elliptic PDE for Stream} satisfies the following estimate:
\begin{equation}\label{Elliptic Estimate for u}
\lVert \psi\rVert_{H^4}\lesssim P(\lVert\rho\rVert_{H^3}) \lVert\omega\rVert_{H^2}.   
\end{equation}
Here $P$ is some monotonically increasing polynomial function.
\end{lem} 

We now state our local existence result.
\begin{thm}\label{LWP of IIE}
Assume $d=2$, Let $(\omega_0,\rho_0)\in H^{2}\times H^{3}$ and let $\rho$ satisfies the Assumption \ref{Assumption A}, then \eqref{Vorticity IIE} admits a unique solution local in time.
\end{thm} 
\begin{proof}[Proof of Theorem~\ref{LWP of IIE}]
 We consider the following iterative scheme:
\begin{equation}\label{Approximation IIE}
\left\{
\begin{aligned}
&\partial_t\omega^{n+1}+{u^n}\cdot\nabla \omega^{n+1}+\nabla\frac{1}{2}\lvert u^n\rvert^2\times\nabla \rho^{n}=0\\
&\partial_t\rho^{n+1}+u^n\cdot\nabla\rho^{n+1}=0\\
&u^{n+1}=\mathcal K_\rho^{n+1} \omega^{n+1}:=\nabla^\perp\bL_{\rho^{n+1}}^{-1}\omega^{n+1}
\end{aligned}
\right.
\end{equation}
The proof is Lagrangian in nature, we will use Sobolev estimate of transport equation to show existence of the fixed point of system \eqref{Approximation IIE}. The proof is divided into $3$ steps.

\begin{itemize}
    \item []\textbf{Step 1.} We define
\begin{equation*}
 \mathcal U_R=\{(\omega,u,\rho)\in C_tH^{2}\times C_tH_\sigma^{3}\times C_tH^{3}: \lVert\omega\rVert_{C_tH_x^2}\le R_1,\lVert u\rVert_{C_tH_x^3}\le R_2, \lVert\rho\rVert_{C_tH_x^{3}}\le R_3\}  
\end{equation*}
We aim to show that for sufficiently small time $T$ and positive numbers $R_1,R_2,R_3$ that we will choose later, the sequence $\{(\omega^n,u^n,\rho^n)\}$ generated by the iterative scheme \eqref{IIE} stays in $U_R$.

First, notice that:
\begin{equation*}
 \rho^{n+1}(t,\cdot)=\rho_0\circ A_t^n ,\quad \text{where }\dot X_t^n=u_t^n(X_t^n),\quad A_t^n=(X_t^n)^{-1}  
\end{equation*}
and $1/\rho^{n+1}$ is also transported. Hence, we have
\begin{equation*}
\lVert\rho_t^{n+1}\rVert_{H^{3}}, \lVert 1/\rho_t^{n+1}\rVert_{H^{3}}\le C(\rho_0) \exp\bigg(C\int_0^t \lVert u^n\rVert_{H^3} d\tau\bigg)\le C(\rho_0)\exp(CtR_2).
\end{equation*}
Here, $C(\rho_0)$ can be chosen as $C(\rho_0)=\lVert\rho_0\rVert_{H^3}+\lVert\rho_0^{-1}\rVert_{H^3}$. 

Now we consider the estimate for vorticity equation. By Duhamel's principle, we can write down the following Lagrangian solution formula:
\begin{equation*}
 \omega^{n+1}(t,\cdot)=\omega_0\circ A_t^n-\bigg[\int_0^t (\nabla^*u^n u^n\cdot\nabla^\perp\rho^n)(X_\tau^n)d\tau \bigg]\circ A_t^n  
\end{equation*}
\begin{equation*}
   =\omega_0\circ A_t^n-\bigg[\nabla^\perp\rho_0\int_0^t \nabla^*X^n (\nabla^*u^n u^n)(X_\tau^n)d\tau \bigg]\circ A_t^n 
\end{equation*}
Hence, we make use of the algebra structure of $H^2$ and standard result in propagation of Sobolev regularity along smooth flow (for instance, see \cite{BCD11}\cite{MP25}) to estimate:
\begin{equation*}
 \lVert\omega^{n+1}\rVert_{H^2}\le \lVert\omega_0\circ A_t\rVert_{H^2}+  \lVert\nabla \rho_0\circ A_t\rVert_{H^2}\exp(CtR_2) \int_0^t  \lVert\nabla X^n\rVert_{H^2} \lVert(\nabla ^* uu)(X^n)\rVert_{H^2} d\tau 
\end{equation*}

\begin{equation*}
\le\lVert\omega_0\rVert_{H^2}\exp(CtR_2)+R_3\exp(2CtR_2)\int_0^t \exp(2C\tau R_2) R_2^2 dt
\end{equation*}
\begin{equation}
 \le \exp(CtR_2)(\lVert\omega_0\rVert_{H^2}+CR_2^2R_3t)   
\end{equation}

Here we temporarily choose $T$ to be sufficiently small such that 
\begin{equation*}
 {CR_2T}\le 1,\quad  CR_2^2R_3T\le \lVert\omega_0\rVert_{H^2}    
\end{equation*}
where $C>1$ is some generic large number. We leave $R_2$ to be fixed later. Thanks to this choice, for all $t\in [0,T]$ we have:
\begin{equation*}
\lVert \rho^{n+1}\rVert_{H^3}, \lVert 1/\rho^{n+1}\rVert_{H^3}\le eC(\rho_0),\quad \lVert\omega_t^{n+1}\rVert_{H^2}\le 2e\lVert\omega_0\rVert_{H^2}
\end{equation*}
Hence, we can choose $R_1= 2 e\lVert\omega_0\rVert_{H^2}$.

Now we estimate $u^{n+1}$. Given $\lVert \omega^{n+1}\rVert\le 2e\lVert\omega_0\rVert_{H^2}$ and $\lVert \rho^{n+1}\rVert_{H^3}\le eC(\rho_0)$, we have:
\begin{equation*}
\lVert u^{n+1}\rVert_{H^3}\lesssim \lVert \bL_{\rho^{n+1}}^{-1}\omega^{n+1}\rVert_{H^4}=\lVert\psi^{n+1}\rVert_{H^4}\le 2 eC P(eC(\rho_0))\lVert\omega_0\rVert_{H^2}    
\end{equation*}
where $P$ is the same polynomial as in \eqref{Elliptic Estimate for u}. Hence, we may choose:
\begin{equation*}
 R_2= 2 eC P( eC(\rho_0))\lVert\omega_0\rVert_{H^2},\quad R_3=eC(\rho_0).  
\end{equation*}
And the above analysis demonstrates that our iteration scheme stays in $\mathcal U_R$.

\item[] \textbf{Step 2.} Now, we show that the above iterative sequence $\{(\omega^n,u^n,\rho^n)\}\subset \mathcal U_R$ is Cauchy in $C_tH^1\times C_tH_\sigma^2\times C_t H^2$. We start by noticing for $u^n=\nabla^\perp \psi^n$:
\begin{equation*}
 \frac{\bL_{\rho^{n+1}}\psi^{n+1}}{\rho^{n+1}}-\frac{\bL_{\rho^{n}}\psi^n}{\rho^n} =\frac{\omega^{n+1}}{\rho^{n+1}}- \frac{\omega^n}{\rho^n}  
\end{equation*}
Which implies:
\begin{equation*}
\Delta (\psi^{n+1}-\psi^n)+\nabla\log\rho^{n+1}\cdot\nabla(\psi^{n+1}-\psi^n)+\nabla(\log\rho^{n+1}-\log\rho^n)\nabla\psi^n=
\end{equation*}
\begin{equation*}
 =\frac{(\omega^{n+1}-\omega^n)}{\rho^{n+1}}+(\frac{1}{\rho^{n+1}}- \frac{1}{\rho^n})\omega^n.    
\end{equation*}
Denote by $C_R$ constants depending on $R_1,R_2, R_3$ which may change from line to line. We now recall the following paraproduct inequality (see \cite{Wi22}): For $f\in H^s$, $g\in H^r$ with $r>s\ge 1$:
\begin{equation}\label{Paraproduct Ineq}
 \lVert f g\rVert_{H^s}\le \lVert f\rVert_{H^s} \lVert g\rVert_{H^r}.   
\end{equation}

By \eqref{Paraproduct Ineq}, we control the RHS as follows:
\begin{equation*}
 \bigg\lVert\frac{(\omega^{n+1}-\omega^n)}{\rho^{n+1}}\bigg\rVert_{H^1}\le C_R \lVert\omega^{n+1}-\omega^n\rVert_{H^1} ,\quad \bigg\lVert(\frac{1}{\rho^{n+1}}- \frac{1}{\rho^n})\omega^n \bigg\rVert_{H^1}\le C_R\lVert \rho^{n+1}-\rho^n\rVert_{H^2},
\end{equation*}
and the commutator term satisfies
\begin{equation*}
 \lVert \nabla(\log\rho^{n+1}-\log\rho^n)\nabla\psi^n\rVert_{H^1}\le C_R\lVert \rho^{n+1}-\rho^n\rVert_{H^2}.
\end{equation*}
Hence, by standard elliptic estimate as in Lemma \ref{Ellptic Regularity}, let $\delta \omega^n:=\omega^{n+1}-\omega^n$, $\delta \rho^n:=\rho^{n+1}-\rho^n$ and $\delta u^n=u^{n+1}-u^n$, we conclude
\begin{equation}\label{delta u n estimate}
 \lVert \delta u^{n}\rVert_{H^2}\le \lVert\psi^{n+1}-\psi^n\rVert_{H^3}\le C_R (\lVert \delta\omega^n\rVert_{H^1}+\lVert \delta\rho^n\rVert_{H^2}) .  
\end{equation}
Straightforward computation yields the following equation for $\delta\omega^n$ and  $\delta \rho^n$:
\begin{equation*}
 \partial_t\delta\omega^n+u^n\cdot\nabla\delta\omega^n+\delta u^n\cdot\nabla\omega^n=\frac{1}{2}(\{\lvert u^n\rvert^2,\rho^n\}-\{\lvert u^{n-1}\rvert^2,\rho^{n-1}\}) =:\sigma_n  ,
\end{equation*}
\begin{equation*}
 \partial_t\delta\rho^n+u^n\cdot\nabla\delta\rho^n+\delta u^n\cdot\nabla\rho^n=0 .  
\end{equation*}
Again, using Duhamel's formula we conclude:
\begin{equation*}
 \delta\omega^n= \delta\omega_0\circ A_t^n+\bigg(\int_0^t (\sigma_n-\delta u^n\cdot\nabla\omega^n)(X_\tau^n)d\tau\bigg)\circ A_t^n ,   
\end{equation*}
\begin{equation*}
 \delta\rho^n= \delta\rho_0\circ A_t^n-\bigg(\int_0^t (\delta u^n\cdot\nabla\rho^n)(X_\tau^n)d\tau\bigg)\circ A_t^n  . 
\end{equation*}
Follow the same line of computation as for $\omega^n$, since $\delta \omega_0=\delta\rho_0=0$, we have the following Sobolev estimates:
\begin{equation*}
 \lVert\delta\omega^n\rVert_{H^1}\le  e\int_0^t \lVert\sigma_n-\delta u^n\cdot\nabla\omega^n\rVert_{H^1}d\tau\le e\int_0^t \lVert\sigma_n\rVert_{H^1}d\tau+C_R\int_0^t \lVert\delta\omega^n\rVert_{H^1}+\lVert\delta \rho^n\rVert_{H^2}d\tau ; 
\end{equation*}
\begin{equation}\label{delta rho n estimate}
 \lVert\delta\rho^n\rVert_{H^2}\le  e\int_0^t \lVert\delta u^n\cdot\nabla\rho^n\rVert_{H^2}d\tau \le C_R\int_0^t \lVert\delta\omega^n\rVert_{H^1}+\lVert\delta \rho^n\rVert_{H^2}d\tau  .
\end{equation}
Now we estimate $\sigma_n$. Notice that
\begin{equation*}
 \{\lvert u^n\rvert^2,\rho^n\}-\{\lvert u^{n-1}\rvert^2,\rho^{n-1}\}=\{\lvert u^n\rvert^2,\rho^n\}-\{\lvert u^{n}\rvert^2,\rho^{n-1}\}+\{\lvert u^{n}\rvert^2,\rho^{n-1}\}-\{\lvert u^{n-1}\rvert^2,\rho^{n-1}\}  
\end{equation*}
\begin{equation*}
=\{\lvert u^n\rvert^2,\delta\rho^{n-1}\}+\nabla^*u^n(u^n-u^{n-1})\nabla^\perp\rho^{n-1}+\nabla^*(u^n-u^{n-1})u^{n-1}\nabla^\perp \rho^{n-1}   \end{equation*}
\begin{equation*}
= \{\lvert u^n\rvert^2,\delta\rho^{n-1}\}+ \nabla^*u^n\delta u^{n-1}\nabla^\perp\rho^{n-1}+ \nabla^*\delta u^{n-1}u^{n-1}\nabla^\perp \rho^{n-1} .
\end{equation*}
Hence, we conclude that
\begin{equation*}
 \lVert\sigma_n\rVert_{H^1}\le C_R (\lVert\delta\omega^{n-1}\rVert_{H^1}+\lVert\delta \rho^{n-1}\rVert_{H^2})   
\end{equation*}
and
\begin{equation}\label{delta omega n estimate}
 \lVert\delta\omega^n\rVert_{H^1}\le C_R \int_0^t \lVert\delta\omega^{n-1}\rVert_{H^1}+\lVert\delta \rho^{n-1}\rVert_{H^2} d\tau+C_R\int_0^t \lVert\delta\omega^{n}\rVert_{H^1}+\lVert\delta \rho^{n}\rVert_{H^2} d\tau.
\end{equation}

\item []\textbf{Step 3.} Now we define: $E^n(t):=\lVert\delta\omega^n\rVert_{H^1}+\lVert\delta\rho^{n}\rVert_{H^2}$. By \eqref{delta u n estimate}, \eqref{delta rho n estimate}, \eqref{delta omega n estimate}, we obtain:
\begin{equation*}
   E^n(t)\le C_R \int_0^t E^n(\tau)d\tau+C_R \int_0^t E^{n-1}(\tau)d\tau,\quad E(0)=0   
\end{equation*}
Now we denote by $U_n:=\sup_{t\in [0,T]} E^n(t)$ and rewrite the above inequality as 
\begin{equation*}
 E^n(t)\le C_R\int_0^t  E^n(\tau)d\tau+C_R \int_0^t E^{n-1}(\tau)d\tau\le C_RTU_n+C_RTU_{n-1}  .
\end{equation*}
Therefore, taking supremum over $t$ on the LHS we conclude:
\begin{equation*}
 U_n\le   C_RTU_n+C_RTU_{n-1} .
\end{equation*}
Choosing $T$ smaller if necessary, such that $C_RT\le 1/3$, we have
\begin{equation*}
 \frac{2}{3}U_n\le  (1-C_RT)U_n\le  C_RTU_{n-1} \le \frac{1}{3}U_{n-1}.
\end{equation*}
Therefore, we obtain $U_n\le U_{n-1}/2\le 2^{-(n-1)}U_1 $. In particular, we demonstrate the summability of the sequence 
\begin{equation*}
 U_n=\sup_{t\in [0,T]} \lVert\delta\omega^n\rVert_{H^1}+\lVert\delta\rho^{n}\rVert_{H^2}.   
\end{equation*}

As a matter of fact, we conclude that $(\omega^n,\rho^n)$ is Cauchy in $C([0,T];H^1\times H^2)$. Now, $u^n$ is also Cauchy in $C([0,T];H_\sigma^2)$, thanks to \eqref{delta u n estimate}. We therefore have a unique limit $(\omega,u,\rho)$.

Finally, by our construction and thanks to Banach-Alaoglu theorem, the sequence $(\omega^n,\rho^n,u^n)$ converges in weak-* sense in $\mathcal U_R$, whose limit must coincide with the strong $C([0,T];H^1\times H^2\times H ^2)$-limit $(\omega,u,\rho)\in \mathcal U_R$, which is the unique solution of \eqref{Vorticity IIE}. 
\end{itemize}
\end{proof}

\section{Lagrangian Analyticity of Inhomogeneous Incompressible Fluids}

 Lagrangian trajectories of many incompressible hydrodynamic models are known to be locally analytic in time, as long as the initial data lies in the local well-posedness class. We refer the reader to \cite{CVW15}\cite{GlassSueurTakahashi2012}\cite{Nadirashvili2013}\cite{Shnirelman2012}\cite{[FZ14]} for detailed treatment on this subject. On the other hand, such a nice property fails to hold when we switch to the compressible models (see \cite{He19}), since the  pressure is no longer an elliptic constraint, but forces the system to be hyperbolic conservation law with finite speed of propagation. Therefore, there may exist points which stay rest and then start to move, which prevents analyticity of the Lagrangian displacement. However, to the author's knowledge, as an intermediate model between incompressible Euler and compressible Euler, whether such a property is enjoyed by IIE is not known in existing literature.

In this section, we approach the Lagrangian analyticity of the inhomogeneous incompressible Euler equation using our Lagrangian formulation. Our approach is similar to \cite{[FZ14]}, where the vorticity formulation is used and Cauchy's invariant plays a key role towards the iterative scheme for Taylor coefficients. In the present work, we stick to the momentum level, which  allows us to deal with general dimensions in a unified way. The main theorem of this section is as follows.

\begin{thm}\label{Lagrangian Analyticity of IIE}
Let $s>d/2$ and assume $u_0\in H_\sigma^{s+1}$, $\rho_0\in H^{s+1}$ satisfies assumption (A), then there exists a time threshold $T_*$ such that the Lagrangian flow $X$ of solution $u$ of IIE is analytic for $t\in (0,T_*)$
\end{thm} 

It turns out that inhomogeneity doesn't prevent analyticity, by comparing above theorem \ref{Lagrangian Analyticity of IIE} with the compressible Euler counterexample in \cite{He19}. The crucial structure of incompressible fluids that leads to the analyticity is the hidden ellipticity in these models (The Biot-Savart law). Such a regularization mechanism, from the Lagrangian viewpoint, is manifested by the following pointwise integration by part lemma (cf. \cite{Co01a}\cite{CI08}).
\begin{lem}\label{Pointwise IBP}
Let $\ell,u\in H^{s+1}$, $s>d/2$. Then $\bP(\nabla^*\ell u)\in H^{s+1}$.
\end{lem} 
\begin{proof}[Proof of Lemma~\ref{Pointwise IBP}]
 It suffices to show $\partial_i\bP(\nabla^*\ell u)\in H^s$ for all $i$. For each fixed $1\le i\le d$, by commutativity between Fourier multipliers we have
\begin{equation*}
 \partial_i \bP(\partial_k\ell_ju_j)=\bP( \partial_{ki}\ell u_j+\partial_k\ell_j\partial_i u_j)=\bP(\partial_k(\partial_i\ell_ju_j)-\partial_i\ell_j\partial_k u_j+\partial_k\ell_j\partial_i u_j) 
\end{equation*}
\begin{equation*}
 =\bP (-\partial_i\ell_j\partial_k u_j+\partial_k\ell_j\partial_i u_j)\in H^s 
\end{equation*}
where the last line follows from the fact that $H^s$ is an algebra and $\bP$ is a $0$-order multiplier. 
\end{proof}

We now proceed to prove theorem \ref{Lagrangian Analyticity of IIE}.

\begin{proof}
We will use our Lagrangian representation to obtain a recursive Taylor formula for the formal time series of Lagrangian displacement, and prove convergence of the formal time series. The proof is divided into $4$ steps.

\begin{itemize}
    \item []
\textbf{Step 1.(Set up the Lagrangian Framework)}

Recall that in proof of theorem 3.1 (III) we derived:
\begin{equation*}
\eta_t=\rho_t u_t+\nabla q_t=\nabla^*A_t\bigg(\eta_0-\int_0^t \frac{1}{2}\lvert \dot X_\tau(x)\rvert^2\nabla\rho_0(x)d\tau\bigg)\circ A_t  
\end{equation*}
Acting $X_t$ to the above equality in the sense of pull-back of covector we have
\begin{equation*}
 \nabla^*X_t(\rho_t u_t+\nabla q_t)\circ X_t= \eta_0-\nabla \rho_0\int_0^t \frac{1}{2}\lvert \dot X_\tau\rvert^2 d\tau
\end{equation*}
Denote by $\ell_t(a):=X_t(a)-a$ the Lagrangian displacement, simplify the LHS, we end up with:
\begin{equation}
\rho_0(\nabla^*\ell+\bId)\dot\ell+\nabla(q_t\circ X_t)+\frac{\nabla\rho_0}{2}\int_0^t \lvert\dot\ell\rvert(\tau)d\tau=\eta_0  ,   
\end{equation}
which suggests that
\begin{equation}\label{Lagrangian Displacement Formula}
\bP\bigg[\rho_0(\nabla^*\ell+\bId)\dot\ell+\frac{\nabla\rho_0}{2}\int_0^t \lvert\dot\ell\rvert^2(\tau)d\tau\bigg]=\bP\eta_0.
\end{equation}
Meanwhile, the following determinant identity relation automatically holds:
\begin{equation}\label{Determinant}
\det(\bId+\nabla\ell(t))=1 \quad\text{for all }t\ge 0   
\end{equation}
Now, assuming for sufficiently small $t>0$, we have the following formal series expansion  of $\ell$:
\begin{equation*}
 \ell_j(t)=\sum_{k=1}^\infty t^k\ell_j^{(k)},   
\end{equation*}
here $\ell_j$ is the $j$-th component of $\ell$. Notice that here $\ell^{(0)}=0$, $\ell^{(1)}=u_0$. Straightforward computation yields:
\begin{equation*}
(\nabla^*\ell(t))_{ij}=\sum_{k=1}^\infty t^k\partial_i\ell_j^{(k)} ,\quad \dot\ell_j(t)=\sum_{k=0}^\infty (k+1)t^k\ell_j^{(k+1)} .
\end{equation*}
Hence, we moreover have:
\begin{equation*}
 (\nabla^*\ell\dot \ell)_i= \sum_{m=1}^\infty t^m\sum_{k=1}^m  (m-k+1)\partial_i\ell_j^{(k)}\ell_j^{(m+1-k)};
\end{equation*}
\begin{equation*}
 \int_0^t\lvert\dot\ell(\tau)\rvert^2d\tau=\sum_{m=1}^\infty t^m\sum_{k=1}^m\frac{k(m+1-k)}{m}\ell_j^{(m+1-k)}\ell_j^{(k)}.
\end{equation*}
Combining the above together, we conclude that \eqref{Lagrangian Displacement Formula} reads for all $n\ge 1$:
\begin{equation*}
 \bP(\rho_0 u_0)=   
\end{equation*}
\begin{equation*}
=\bP \bigg[\rho_0\ell_i^{(1)}+\sum_{m=1}^\infty t^m\bigg((m+1)\ell_i^{(m+1)}\rho_0+\sum_{k=1}^m (m+1-k)\big(\partial_i\ell_j^{(k)}\ell_j^{(m-1+k)}\rho_0+\frac{k\partial_i \rho_0}{2m}\ell_j^{(m+1-k)}\ell_j^{(k)}\big)\bigg)\bigg],
\end{equation*}
which implies that coefficient of $t^m$ vanishes at all order: For all $m\ge 1$
\begin{equation}\label{Taylor Coefficient 1}
 \bP\bigg[(m+1)\ell_i^{(m+1)}\rho_0+\sum_{k=1}^m (m+1-k)\big(\partial_i\ell_j^{(k)}\ell_j^{(m-1+k)}\rho_0+\frac{k\partial_i \rho_0}{2m}\ell_j^{(m+1-k)}\ell_j^{(k)}\big)\bigg]=0  .
\end{equation}
Moreover, the determinant identity \eqref{Determinant} naturally implies
\begin{equation*}
0=\log\det(\bId+\nabla \ell(t))=\Tr\log(\bId+\nabla\ell(t))=\Tr\bigg(\sum_{k=1}^\infty \frac{(-1)^k}{k} [\nabla \ell(t)]^k \bigg) 
\end{equation*}
\begin{equation*}
= \sum_{k=1}^\infty \frac{(-1)^k}{k} \Tr\bigg[\bigg(\sum_{r=1}^\infty \nabla\ell^{(r)}t^r\bigg)^k\bigg] .
\end{equation*}
Hence, the vanishing of coefficient for $t^m$ implies:
\begin{equation}\label{Taylor Coefficient 2}
 0=\sum_{k=1}^m\frac{(-1)^k}{k}\sum_{r_1+\cdots+r_k=m, r_i\ge 1}\Tr(\nabla \ell^{(r_1)}\cdots\nabla\ell^{(r_k)}).
\end{equation}
\item[] \textbf{Step 2.} We now derive the iterative scheme for coefficients $\{\ell^{(k)}\}_{k\ge 1}$ based on \eqref{Taylor Coefficient 1}\eqref{Taylor Coefficient 2}. To this end, we rewrite \eqref{Taylor Coefficient 1} as
\begin{equation}\label{Recursive l}
 (m+1)\ell_i^{(m+1)}\rho_0+U_i^{(m)}=\nabla p_i^{(m)}  
\end{equation}
where $U^{(m)}$ is a vector field quadratic in lower-order terms:
\begin{equation*}
  U^{(m)}= \sum_{k=1}^m (m+1-k)\big(\partial_i\ell_j^{(k)}\ell_j^{(m-1+k)}\rho_0+\frac{k\partial_i \rho_0}{2m}\ell_j^{(m+1-k)}\ell_j^{(k)}\big) .
\end{equation*}
Meanwhile, \eqref{Taylor Coefficient 2} implies:
\begin{equation}
\Tr(\nabla\ell^{(m)})=\nabla\cdot\ell^{(m)}=\sum_{k=2}^m\frac{(-1)^k}{k}\sum_{r_1+\cdots+r_k=m,r_i\ge 1}\Tr(\nabla \ell^{(r_1)}\cdots\nabla\ell^{(r_k)})=:W^{(m)}   
\end{equation}
Here, the expansion eventually leads to a polynomial of lower order $\nabla\ell^{(r)}$s' entries with degree at most $m$.

Therefore, the following recursive relation holds: Fix $m$ and assume we already know $\{\ell^{(k)}\}_{k=1}^m$, then we can explicitly compute $U^{(m)}$, $W^{(m+1)}$ and $\nabla\cdot \ell^{(m+1)}=W^{(m+1)}$. Taking into account \eqref{Recursive l} we have:
\begin{equation*}
  \nabla\cdot\ell^{(m+1)}=\frac{1}{m+1}\nabla\cdot(\rho_0^{-1}(\nabla p^{(m)}-U^{(m)}))=W^{(m+1)}  ,
\end{equation*}
which guarantees:
\begin{equation}\label{Pressure}
  \nabla p^{(m)}= \nabla \Delta_{\rho_0}^{-1}[(m+1)W^{(m+1)}+\nabla\cdot(\rho_0^{-1}U^{(m)})].
\end{equation}
Therefore one can produce $\ell^{(m+1)}$ from $\{\ell^{(k)}\}_{k=1}^m$ as follows:
\begin{equation*}
\ell^{(m+1)}=  \frac{1}{(m+1)\rho_0}\big[\nabla \Delta_{\rho_0}^{-1}[(m+1)W^{(m+1)}-\nabla\cdot(\rho_0^{-1}U^{(m)})]-U^{(m)}\big] 
\end{equation*}
Which is an explicit recursion formula for Taylor coefficients.

\item[] \textbf{Step 3.(Elliptic Estimates)} Introduce the following generating function:
\begin{equation}\label{Series of l}
\bar\zeta(t):=\sum_{k=1}^\infty t^k \lVert\nabla \ell^{(k)}\rVert_{H^s}    
\end{equation}
We aim to show that the series \eqref{Series of l} converges absolutely.

Here, we consider the Leray-Hodge decomposition $U^{(m)}=\bP U^{(m)}+\nabla\varphi^{(m)}$. Notice then
\begin{equation*}
\nabla\Delta_{\rho_0}^{-1}\nabla\cdot(\rho_0^{-1}U^{(m)})=\nabla\varphi^{(m)}+ \nabla\Delta_{\rho_0}^{-1}\nabla\cdot(\rho_0^{-1}\bP U^{(m)})   
\end{equation*}

First, we need to show that $\ell^{(m+1)}\in H^{s+1}$, given $\{\ell^{(k)}\}_{k=1}^m\subset H^{s+1}$. Note that with our assumption \ref{Assumption A}, $\Delta_{\rho_0}$ is uniformly elliptic. The above equality \eqref{Pressure} reads
\begin{equation}\label{Pressure*}
  \nabla (p^{(m)}-\varphi^{(m)})= \nabla \Delta_{\rho_0}^{-1}[(m+1)W^{(m+1)}-\nabla\cdot(\rho_0^{-1}\bP U^{(m)})]
\end{equation}
and we have by \eqref{Recursive l}:
\begin{equation}\label{l^m Reformulate}
\ell^{(m+1)}=\frac{1}{(m+1)\rho_0}\big[-\bP U^{(m)}+\nabla(p^{(m)}-\varphi^{(m)})\big]    .
\end{equation}

We first control $\bP U^{(m)}$. We notice that:
\begin{equation*}
 \bP U^{(m)}=\sum_{k=1}^m(m-1+k)\big(\bP[\nabla^*\ell^{(k)}(\ell^{(m+1-k)}\rho_0)]+\frac{k}{2m}\bP[\nabla\rho_0(\ell^{(m+1-k)}\cdot\ell^{(k)}) ])  
\end{equation*}

Therefore, $\bP U^{(m)}\in H^{s+1}$ will hold true if we succeed in showing for any $1\le k\le m$:
\begin{equation*}
\bP[\nabla^*\ell^{(k)}(\ell^{(m+1-k)}\rho_0)]\in H^{s+1},\quad   \bP[\nabla\rho_0(\ell^{(m+1-k)}\cdot\ell^{(k)}) ]\in H^{s+1}.  
\end{equation*}
Combining Lemma 4.2 and a simple manipulation on the second term:
\begin{equation*}
\bP[\nabla\rho_0(\ell^{(m+1-k)}\cdot\ell^{(k)}) ]=-\bP[\nabla(\ell^{(m+1-k)}\cdot\ell^{(k)})\rho_0]   \end{equation*}
\begin{equation*}
 =-\bP[\nabla^*\ell^{(m+1-k)}\cdot\ell^{(k)}\rho_0]-\bP[\nabla^*\ell^{(k)}\cdot\ell^{(m+1-k)}\rho_0],
\end{equation*}
we eventually conclude the following estimate:
\begin{equation}\label{PU Estimate}
 \lVert\bP U^{(m)}  \rVert_{H^{s+1}}\lesssim\sum_{k+j=m+1; k,j\ge 1} \lVert\nabla \ell^{(k)} \rVert_{H^{s}} \lVert\nabla \ell^{(j)} \rVert_{H^{s}}.
\end{equation}

On the other hand, rewrite \eqref{Pressure} as
\begin{equation*}
\nabla (p^{(m)}-\varphi^{(m)})=\nabla\Delta_{\rho_0}^{-1}[(m+1)W^{(m+1)}-\nabla\cdot(\rho_0^{-1}\bP U^{(m)})] .   
\end{equation*}
Now, a slight variant of Lemma \ref{Ellptic Regularity} (taking $\rho=1/\rho_0$) implies:
\begin{equation*}
\lVert\nabla (p^{(m)}-\varphi^{(m)})\rVert_{H^{s+1}}\le C(\rho_0)(m+1)\lVert W^{(m+1)}\rVert_{H^{s}}+ \lVert   \nabla\Delta_{\rho_0}^{-1} \nabla\cdot(\rho_0^{-1}\bP U^{(m)})\rVert_{H^{s+1}}.
\end{equation*}
Since $W^{(m+1)}$ consists of the sum of product of $\nabla\ell^{(r)}$ with $1\le r\le m$, the fact that $H^s$ is an algebra implies $W^{(m+1)}\in H^s$, with: 
\begin{equation}
 \lVert W^{(m+1)}\rVert_{H^s}\lesssim \sum_{j=2}^d\sum_{r_1+\cdots+r_j=m+1, r_i\ge 1} \prod_{i=1}^j\lVert\nabla \ell^{(r_i)}\rVert_{H^s},
\end{equation}
while above estimate \eqref{PU Estimate} implies
\begin{equation*}
\lVert \nabla\Delta_{\rho_0}^{-1}\nabla\cdot(\rho_0^{-1}\bP U^{(m)})\rVert_{H^{s+1}}\lesssim \sum_{k+j=m+1; k,j\ge 1} \lVert\nabla \ell^{(k)} \rVert_{H^{s}} \lVert\nabla \ell^{(j)} \rVert_{H^{s}} .   
\end{equation*}

Hence, we conclude from \eqref{l^m Reformulate}:
\begin{equation}
\lVert\nabla \ell^{(m+1)}\rVert_{H^s}\lesssim  \sum_{k+j=m+1; k,j\ge 1} \lVert\nabla \ell^{(k)} \rVert_{H^{s}} \lVert\nabla \ell^{(j)} \rVert_{H^{s}} + \sum_{j=2}^d\sum_{r_1+\cdots+r_j=m+1} \prod_{i=1}^j\lVert\nabla \ell^{(r_i)}\rVert_{H^s}  , 
\end{equation}
where RHS is a polynomial of $\lVert\nabla\ell^{(r)}\rVert_{H^s}$ up to order $d$. Hence, there exists positive numbers $C_2,\cdots,C_d$ such that 
\begin{equation}\label{Controlled Recurrence}
\lVert\nabla \ell^{(m+1)}\rVert_{H^s}\le \sum_{j=2}^d C_j\sum_{r_1+\cdots+r_j=m+1} \prod_{i=1}^j\lVert\nabla \ell^{(r_i)}\rVert_{H^s}    .
\end{equation}
\item[] \textbf{Step 4. (Conclude Analyticity)}
Now based on the recurrence relation, we can recursively build a sequence $\{\beta_m\}$ dominating coefficients $\{\lVert \nabla\ell^{(m)}\rVert_{H^s}\}$ of our generating function $\zeta$ as follows:
\begin{equation*}
 \beta_1=\lVert \nabla u_0\rVert_{H^s},\quad \beta_{m+1}=  \sum_{j=2}^d C_j\sum_{r_1+\cdots+r_j=m+1} \prod_{i=1}^j\beta_{r_i}   
\end{equation*}
and we define the formal series
\begin{equation*}
 \tilde\zeta(t)=\sum_{m=1}^\infty \beta_m t^m   
\end{equation*}
which is the formal expansion of the unique root of the following polynomial:
\begin{equation}\label{Zeta equality}
P(t,\zeta)=-\zeta+t\lVert\nabla u_0\rVert_{H^s}+C_2\zeta^2+\cdots+C_d\zeta^d    
\end{equation}

To see the uniqueness of root of $P(t,\cdot)$, notice first that for each $t\ge 0$, solving $P(t,\zeta)=0$ is equivalent to solve 
\begin{equation*}
J_t(\alpha):=\alpha-\sum_{j=2}^d C_j\alpha^j=t\Vert\nabla u_0\rVert_{H^s}    
\end{equation*}
We observe that for any $t>0$, $J_t(0)=0$, while $J_t^\prime(0)=1>0$ and $J_t(\alpha)\to -\infty$ as $\alpha\to\infty$. Hence, we consider the first  $\alpha_*$ with sufficiently small derivative: Fix $0<\delta\ll1$, we consider
\begin{equation*}
 \alpha^*=\min\{\alpha: J^\prime(\alpha)=\delta\},\quad J^\prime(\alpha)>\delta\text{ for all }0<\alpha<\alpha^*   
\end{equation*}
We now define 
\begin{equation*}
\Psi_t(\zeta):=t\rVert\nabla u_0\rVert_{H^s}+\sum_{j=2}^d C_j\zeta^j 
\end{equation*}
Notice that
\begin{equation*}
\lvert \Psi_t^\prime(\zeta)\rvert\le \sum_{j=2}^d jC_j \zeta^{j-1} = 1-J^\prime(\zeta)\le 1-\delta,\quad\text{for all }\zeta\in [0,\alpha_*]
\end{equation*}
Hence, for each $t$, $\Psi_t(\zeta)$ is a contraction. Moreover, we notice that $\Psi_t$ is increasing, therefore:
\begin{equation*}
\Psi_t([0,\alpha])\subset [\Psi_t(0),\Psi_t(\alpha)]=\bigg[t\Vert \nabla u_0\rVert_{H^s},t\Vert \nabla u_0\rVert_{H^s}+\sum_{j=2}^d C_j\alpha^j \bigg]    
\end{equation*}
To apply Banach fixed point theorem, it suffices to show:
\begin{equation*}
 t\Vert \nabla u_0\rVert_{H^s}+\sum_{j=2}^d C_j\alpha^j  \le \alpha 
\end{equation*}
Hence it suffices to pick $t$ such that
\begin{equation*}
  t\le T(\alpha):=\Vert\nabla u_0\rVert_{H^s}^{-1} J(\alpha)
\end{equation*}
Therefore, for all $t\in [0,T(\alpha)]$ with $\alpha\in (0,\alpha^*)$, $\Psi_t$ has a unique fixed point. Hence, we find a unique $\zeta_*\in (0,\alpha^*)$ such that $P(t,\zeta^*)=0$. 

Now, since $P(t,\zeta)$ satisfies
\begin{equation*}
 P(0,0)=0,\quad \partial_\zeta P(0,0)=-1\neq 0   
\end{equation*}
The implicit function theorem implies that locally around $(0,0)$, there exists an analytic function $\zeta_*(t)$ such that $P(t,\zeta_*(t))=0$. By uniqueness of root of $P(t,\cdot)$, we conclude the formal series solution must coincide with analytic function $\zeta_*(t)$, which concludes the absolute convergence of series $\tilde\zeta$ in $t\in (0,T(\alpha^*))$, hence also the convergence of $\bar\zeta(t)$. Hence, proof of the Lagrangian analyticity of IIE with convergence radius $t\in (0,T(\alpha^*))$ is complete. 
\end{itemize}
\end{proof}

\bibliographystyle{siam}
\bibliography{refs}

@article{GlassSueurTakahashi2012,
  author  = {Glass, Olivier and Sueur, Franck and Takahashi, Tak{\'e}o},
  title   = {Smoothness of the motion of a rigid body immersed in an incompressible perfect fluid},
  journal = {Annales scientifiques de l'{\'E}cole Normale Sup{\'e}rieure},
  volume  = {45},
  number  = {1},
  pages   = {1--51},
  year    = {2012}
}

@article{Shnirelman2012,
  author  = {Shnirelman, Alexander},
  title   = {On the Analyticity of Particle Trajectories in the Ideal Incompressible Fluid},
  journal = {arXiv preprint arXiv:1205.5837},
  year    = {2012}
}

@article{Nadirashvili2013,
  author  = {Nadirashvili, N.},
  title   = {On stationary solutions of two-dimensional Euler Equation},
  journal = {Archive for Rational Mechanics and Analysis},
  volume  = {209},
  number  = {3},
  pages   = {729--745},
  year    = {2013}
}

@article{Ar66,
  author  = {Arnold, Vladimir},
  title   = {Sur la g{\'e}om{\'e}trie diff{\'e}rentielle des groupes de {L}ie de dimension infinie et ses applications {\`a} l’hydrodynamique des fluides parfaits},
  journal = {Annales de l'Institut Fourier},
  volume  = {16},
  number  = {1},
  pages   = {319--361},
  year    = {1966}
}

@book{AK98,
  author    = {Arnold, Vladimir I. and Khesin, Boris A.},
  title     = {Topological Methods in Hydrodynamics},
  publisher = {Springer},
  address   = {New York},
  year      = {1998}
}

@book{BCD11,
  author    = {Bahouri, Hajer and Chemin, Jean-Yves and Danchin, Rapha{\"e}l},
  title     = {Fourier Analysis and Nonlinear Partial Differential Equations},
  series    = {Grundlehren der Mathematischen Wissenschaften},
  volume    = {343},
  publisher = {Springer},
  address   = {Berlin, Heidelberg},
  year      = {2011},
  isbn      = {978-3-642-16829-1},
  doi       = {10.1007/978-3-642-16830-7}
}

@inproceedings{BCHM00,
  author    = {Bloch, Anthony and Crouch, Peter E. and Holm, Darryl D. and Marsden, Jerrold E.},
  title     = {An optimal control formulation for inviscid incompressible ideal fluid flow},
  booktitle = {Proceedings of the 39th IEEE Conference on Decision and Control},
  volume    = {2},
  pages     = {1273--1278},
  year      = {2000},
  publisher = {IEEE}
}

@article{[BLS20],
  title   = {A blow-up criterion for the inhomogeneous incompressible Euler equations},
  author  = {Bae, Hantaek and Lee, Woojae and Shin, Jaeyong},
  journal = {Nonlinear Analysis},
  volume  = {196},
  pages   = {111774},
  year    = {2020},
  doi     = {10.1016/j.na.2020.111774}
}

@article{Br89,
  author  = {Brenier, Yann},
  title   = {The least action principle and the related concept of generalized flows for incompressible perfect fluids},
  journal = {Journal of the American Mathematical Society},
  volume  = {2},
  number  = {2},
  pages   = {225--255},
  year    = {1989}
}

@article{Br99,
  author  = {Brenier, Yann},
  title   = {Minimal geodesics on groups of volume-preserving maps and generalized solutions of the {E}uler equations},
  journal = {Communications on Pure and Applied Mathematics},
  volume  = {52},
  number  = {4},
  pages   = {411--452},
  year    = {1999}
}

@article{Ch03,
  author  = {Chae, Dongho},
  title   = {Local Existence and Blow-up Criterion of the Inhomogeneous {E}uler Equations},
  journal = {Journal of Mathematical Fluid Mechanics},
  volume  = {5},
  number  = {1},
  pages   = {1--18},
  year    = {2003}
}

@article{Co01a,
  author  = {Constantin, Peter},
  title   = {An {E}uler--{L}agrangian approach for incompressible fluids: local theory},
  journal = {Journal of the American Mathematical Society},
  volume  = {14},
  number  = {2},
  pages   = {263--278},
  year    = {2001}
}

@article{CVW15,
title = {Analyticity of Lagrangian trajectories for well posed inviscid incompressible fluid models},
journal = {Advances in Mathematics},
volume = {285},
pages = {352-393},
year = {2015},
issn = {0001-8708},
doi = {https://doi.org/10.1016/j.aim.2015.05.019},
url = {https://www.sciencedirect.com/science/article/pii/S0001870815003011},
author = {Peter Constantin and Vlad Vicol and Jiahong Wu},}

@article{CI08,
  author  = {Constantin, Peter and Iyer, Gautam},
  title   = {A stochastic {L}agrangian representation of the three-dimensional incompressible {N}avier--{S}tokes equations},
  journal = {Communications on Pure and Applied Mathematics},
  volume  = {61},
  number  = {3},
  pages   = {330--345},
  year    = {2008}
}

@misc{ElgindiPasqualotto2023,
  author       = {Tarek M. Elgindi and Federico Pasqualotto},
  title        = {From Instability to Singularity Formation in Incompressible Fluids},
  year         = {2023},
  eprint       = {2310.19780},
  archivePrefix= {arXiv},
  primaryClass = {math.AP},
  url          = {https://arxiv.org/abs/2310.19780}
}

@article{[EP25],
  author  = {Elgindi, Tarek M. and Pasqualotto, Federico},
  title   = {Invertibility of a Linearized Boussinesq Flow: A Symbolic Approach},
  journal = {Communications in Mathematical Physics},
  year    = {2025},
  volume  = {406},
  number  = {11},
  doi     = {10.1007/s00220-025-05367-6}
}

@article{Da06,
  author  = {Danchin, Rapha{\"e}l},
  title   = {The inviscid limit for density-dependent incompressible fluids},
  journal = {Annales de la Facult{\'e} des Sciences de Toulouse. Math{\'e}matiques},
  volume  = {15},
  number  = {4},
  pages   = {637--688},
  year    = {2006}
}

@article{[DF11],
  title   = {The well-posedness issue for the density-dependent Euler equations in endpoint Besov spaces},
  author  = {Danchin, Rapha{\"e}l and Fanelli, Francesco},
  journal = {J. Math. Pures Appl. (9)},
  volume  = {96},
  number  = {3},
  pages   = {253--278},
  year    = {2011},
  doi     = {10.1016/j.matpur.2011.04.005}
}

@article{[Fa25],
  title         = {Geometric blow-up criteria for the non-homogeneous incompressible Euler equations in 2-D},
  author        = {Fanelli, Francesco},
  year          = {2025},
  eprint        = {2502.10024},
  archivePrefix = {arXiv},
  primaryClass  = {math.AP},
  note          = {arXiv:2502.10024}
}

@article{[Da10],
  title   = {On the well-posedness of the incompressible density-dependent Euler equations in the $L^p$ framework},
  author  = {Danchin, Rapha{\"e}l},
  journal = {Journal of Differential Equations},
  volume  = {248},
  number  = {8},
  pages   = {2130--2170},
  year    = {2010},
  doi     = {10.1016/j.jde.2009.09.007}
}

@book{[Lio96],
  title     = {Mathematical Topics in Fluid Mechanics. Volume 1: Incompressible Models},
  author    = {Lions, Pierre-Louis},
  series    = {Oxford Lecture Series in Mathematics and its Applications},
  number    = {3},
  address   = {Oxford},
  publisher = {Clarendon Press, Oxford University Press},
  year      = {1996}
}

@misc{DrWS25,
  author       = {Drivas, Theodore D.},
  title        = {Winter School: Boundary and Singularity in Fluid Mechanics — Drivas Course, Lecture 1},
  year         = {2025},
  howpublished = {\url{https://www.math.stonybrook.edu/~tdrivas/Courses/winterschool25/DrivasLecture1.pdf}},
  note         = {Lecture notes, Simons Center for Geometry and Physics, Stony Brook University; winter school held Jan 6--10, 2025},
  urldate      = {2025-11-03}
}

@article{DG26+,
  author  = {Drivas, Theodore D. and Danill Glukhovskiy},
 title= {On the Realization of Ideal Constraints and
Incompressible/Inextensible Limits with Ill-Prepared Initial
Data},
year = {In Preparation}
}

@article{EM70,
  author  = {Ebin, David G. and Marsden, Jerrold E.},
  title   = {Groups of diffeomorphisms and the motion of an incompressible fluid},
  journal = {Annals of Mathematics},
  volume  = {92},
  number  = {1},
  pages   = {102--163},
  year    = {1970}
}

@article{[FZ14],
  title   = {Time-analyticity of Lagrangian particle trajectories in ideal fluid flow},
  author  = {Zheligovsky, Vladislav and Frisch, Uriel},
  journal = {Journal of Fluid Mechanics},
  volume  = {749},
  pages   = {404--430},
  year    = {2014},
  doi     = {10.1017/jfm.2014.221}
}

@misc{Ho21,
  author       = {Holm, Darryl D.},
  title        = {MATH97178 Lecture Notes: Dynamics, Symmetry and Integrability},
  year         = {2021},
  institution  = {Imperial College London},
  howpublished = {\url{https://www.ma.ic.ac.uk/~dholm/classnotes/GM2-Spring-2021-notesRev1.pdf}},
  note         = {Lecture notes, Spring Term 2021},
  urldate      = {2025-12-01}
}

@article {He19,
author= {Hernandez, Matthew},
title={Mechanisms of Lagrangian Analyticity in Fluids},
Journal={Archive for Rational Mechanics and Analysis},
Year={2019},
Volume={233},
Pages={513--598}
}

@article{[BKM84],
  author  = {Beale, J. T. and Kato, Tosio and Majda, Andrew},
  title   = {Remarks on the breakdown of smooth solutions for the 3-D Euler equations},
  journal = {Communications in Mathematical Physics},
  volume  = {94},
  number  = {1},
  pages   = {61--66},
  year    = {1984},
  doi     = {10.1007/BF01212349}
}

@article{HMR98,
  author  = {Holm, Darryl D. and Marsden, Jerrold E. and Ratiu, Tudor S.},
  title   = {The {E}uler--{P}oincar{\'e} equations and semidirect products with applications to continuum mechanics},
  journal = {Advances in Mathematics},
  volume  = {137},
  pages   = {1--81},
  year    = {1998}
}

@misc{Kh25,
  author  = {Boris Khesin},
 title= {Personal Communication}
}

@article{KMM20,
  author  = {Khesin, Boris and Misio{\l}ek, Gerard and Modin, Klas},
  title   = {Geometric hydrodynamics and infinite-dimensional Newton's equations},
  journal = {Bulletin of the American Mathematical Society},
  volume  = {58},
  pages   = {377--442},
  year    = {2021}
}

@article{[Mar76],
  author  = {Marsden, Jerrold E.},
  title   = {Well-posedness of the equations of a non-homogeneous perfect fluid},
  journal = {Communications in Partial Differential Equations},
  year    = {1976},
  volume  = {1},
  number  = {3},
  pages   = {215--230},
  doi     = {10.1080/03605307608820010}
}

@article{[LLP11],
  author  = {Lopes Filho, Milton C. and Nussenzveig Lopes, Helena J. and Precioso, Juliana C.},
  title   = {Least action principle and the incompressible Euler equations with variable density},
  journal = {Transactions of the American Mathematical Society},
  year    = {2011},
  volume  = {363},
  number  = {5},
  pages   = {2641--2661},
  doi     = {10.1090/S0002-9947-2010-05206-7}
}

@book{MR99,
  author    = {Marsden, Jerrold E. and Ratiu, Tudor S.},
  title     = {Introduction to Mechanics and Symmetry: A Basic Exposition of Classical Mechanical Systems},
  edition   = {2},
  series    = {Texts in Applied Mathematics},
  volume    = {17},
  publisher = {Springer},
  address   = {New York},
  year      = {1999}
}

@article{MW83,
  author  = {Marsden, Jerrold E. and Weinstein, Alan},
  title   = {Coadjoint orbits, vortices, and {C}lebsch variables for incompressible fluids},
  journal = {Physica D: Nonlinear Phenomena},
  volume  = {7},
  number  = {1--3},
  pages   = {305--323},
  year    = {1983}
}

@article{Pa26+,
  author  = {Pan, Anping},
  title   = {Remarks on the Lifespan and Continuation Criteria of Two Dimensional Incompressible Fluid Models.},
 journal = {Preprint, arXiv:2604.15721}
}

@book{TeschlODE2012,
  author    = {Teschl, Gerald},
  title     = {Ordinary Differential Equations and Dynamical Systems},
  series    = {Graduate Studies in Mathematics},
  volume    = {140},
  publisher = {American Mathematical Society},
  address   = {Providence, RI},
  year      = {2012},
  isbn      = {978-0-8218-8328-0}
}

@article{MP25,
  author  = {Mazzucato, Anna and Pan, Anping},
  title   = {Variational Principle and Stochastic Lagrangian Formulation of Hydrodynamic Equations},
 journal = {arXiv:2511.21498}
}

@misc{Wi22,
  author = {Williams, Michael},
  title  = {Notes on Harmonic Analysis},
  year   = {2022},
  note   = {Lecture notes}
}

@article{Zh10,
  author  = {Zhang, Xicheng},
  title   = {A stochastic representation for backward incompressible {N}avier--{S}tokes equations},
  journal = {Probability Theory and Related Fields},
  volume  = {148},
  number  = {1--2},
  pages   = {305--332},
  year    = {2010}
}

\end{document}